\newcommand{\Fix}{\mathrm {Fix}}
\newcommand{\fix}{\mathrm {fix}}
\newcommand{\Per}{\mathrm {Per}}
\newcommand{\trace}{\mathrm {trace\,}}

\def\C{{\mathbb C}}
\def\N{{\mathbb N}}
\def\R{{\mathbb R}}
\def\Z{{\mathbb Z}}

\def\cB{{\mathcal B}}
\def\cL{{\mathcal L}}
\def\cK{{\mathcal K}}
\def\cF{{\mathcal F}}
\def\cG{{\mathcal G}}

\def\id{{\mathrm{id}}}
\newcommand{\oline}{\overline}

\documentclass[a4paper]{article}

\usepackage[latin1]{inputenc}
\usepackage{fontenc}

\font\teneusm=eusm10
\font\seveneusm=eusm7
\font\fiveeusm=eusm5
\newfam\eusmfam
\textfont\eusmfam\teneusm
\scriptfont\eusmfam\seveneusm
\scriptscriptfont\eusmfam\fiveeusm

\usepackage{amsfonts,amsmath,amssymb,theorem}
\usepackage{amsxtra}
\input xy
\xyoption{all}
\accentedsymbol{\Ahatdot}{\Dot{\Hat A}}
\accentedsymbol{\varrhohatdot}{\Dot{\Hat \varrho}}
\accentedsymbol{\Fbarhat}{\Hat{\Bar F}}

\newtheorem{Thm}{Theorem}[section]

\newtheorem{Prop}[Thm]{Proposition}
\newtheorem{Lemma}[Thm]{Lemma}
{\theorembodyfont{\rmfamily} \newtheorem{Rem}[Thm]{Remark} }
{\theorembodyfont{\rmfamily}  }

\def\pf{\noindent {\it Proof.}\hskip 8pt}
\def\qed{\hspace*{\fill}{\setlength{\fboxrule}{0.8pt}\setlength{\fboxsep}{1mm}
\fbox{\null}} \vskip 10pt}

\begin{document}
\title
{Transfer operators and dynamical zeta functions for a class of lattice spin
 models}
\author{J. Hilgert\thanks{
Institut f\"ur Mathematik, Technische Universit\"at Clausthal,
38678 Clausthal-Zellerfeld, Germany.
E-mail: hilgert@math.tu-clausthal.de} \ \ and\
  D. Mayer\thanks{%
IHES, 91440 Bures sur Yvette, France, on leave of absence from Institut
f\"ur Theoretische Physik, Technische Universit\"at Clausthal,
38678 Clausthal-Zellerfeld, Germany. E-mail: mayer@ihes.fr or
dieter.mayer@tu-clausthal.de}}
\date{March 11, 2002}
\maketitle
\makeatother

{\bf Abstract}: 
{\small We investigate the location of zeros and poles of a dynamical zeta function for a family of 
subshifts of finite type with an interaction function depending on the parameters 
$\underline{\lambda}=(\lambda_{1},\ldots,\lambda_{m})$ with $0 \leqslant \lambda_{i} \leqslant 1$.
 The system corresponds to the well known Kac-Baker lattice spin model in statistical mechanics. 
Its dynamical  zeta function can be expressed in terms of the Fredholm determinants of two transfer 
operators $\mathcal{L}_{\beta}$ and $\mathcal{G}_{\beta}$ with $\mathcal{L}_{\beta}$ the Ruelle operator 
acting in a Banach space of holomorphic functions,
 and $\mathcal{G}_{\beta}$ an integral operator introduced originally by Kac, which acts in the space 
$L^{2}(\mathbb{R}^{m}, d\underline{x})$ with a kernel which is symmetric and positive definite for 
positive $\beta$. 
By relating the two operators to each other via the Segal-Bargmann transform  we prove equality of their 
spectra and hence reality, respectively positivity, for the eigenvalues of the operator 
$\mathcal{L}_{\beta}$ for real, respectively positive, $\beta$. 
For a restricted range of parameters  $0 \leqslant \lambda_{i} \leqslant \frac{1}{2}$, 
$1\leqslant i \leqslant m$ we can determine the asymptotic behavior of the eigenvalues of 
$\mathcal{L}_{\beta}$ for large positive and negative values of $\beta$ and deduce from this the existence 
of infinitely many non trivial zeros and poles of the dynamical zeta functions on the real $\beta$ line 
at least for generic $\underline{\lambda}$. 
For the special choice $\lambda_{i}=\frac{1}{2}$, $1\leqslant i \leqslant m$, we find a family of 
eigenfunctions and eigenvalues of $\mathcal{L}_{\beta}$ leading to an infinite sequence of equally 
spaced ``trivial'' zeros and poles of the zeta function on a line parallel to the imaginary $\beta$-axis. 
Hence there seems to hold some generalized Riemann hypothesis also for this kind of dynamical zeta 
functions. 

}

\section{Introduction}

The transfer matrix method has played an important role in statistical mechanics
ever since E.~Ising for the first time used this method to solve his 1-dimensional
lattice spin model with nearest neighbour interaction.
The method was extended later to treat higher dimensional models with arbitrary
finite range interactions.
The most satisfying theory for this method from the mathematical point of view
goes back to D.~Ruelle who introduced the so called transfer operator for 1-dimensional
lattice spin systems with arbitrary long range interactions (see \cite{Ru68}).
Continued interest in such systems is related to the fact that these systems show
up in a rather natural way within the so called thermodynamic formalism
for dynamical systems (see \cite{Ru78}).
Thereby the  transfer operator can be used for instance to construct invariant measures
for such systems and to characterize their ergodic properties (see \cite{Ba00}).

Another nice application of this method is to the theory of dynamical zeta
functions (see \cite{Ru92}, \cite{Ru94}). These functions can be interpreted as generating functions
for the partition functions of the system constructed in complete analogy to the partition
functions of lattice spin systems.
It turned out that the transfer operator  method had indeed been used already some time ago in the
$p$-adic setup of zeta functions by B. Dwork (see \cite{Dw60} or \cite{Ro86}) who constructed such an operator
to show rationality of the Artin-Weil zeta function for projective algebraic
varieties over finite fields and proved in this way part of the Weil conjectures
(see \cite{We49}).
More recently the method also yielded a completely new approach to Selberg's zeta function,
which can also be viewed as a dynamical zeta function for the geodesic
flow on surfaces of constant negative curvature (see \cite{Ma91}).
Indeed, the aforementioned Artin-Weil zeta function is nothing but the dynamical
Artin-Mazur zeta function (see \cite{ArMa65}) for the Frobenius map of the algebraic variety.
Typically, such dynamical zeta functions can be expressed in terms of some kind of Fredholm
determinant of the transfer operator, which therefore allows a spectral interpretation
of the zeros and poles of these functions. The existence of such an interpretation is one
of the challenging open problems for all arithmetic zeta functions
of number theory and algebraic geometry (see \cite{Be86}, \cite{Co96}, \cite{De99}).
Obviously such a spectral interpretation is also closely related to another
famous open problem for these zeta and more general L-functions, namely the general Riemann
hypothesis: one expects that the zeros and poles of such functions are located on
critical lines in the complex plane as one does in the special case of the well known Riemann zeta
function. Presently it is not known whether such a conjecture makes sense also for general
dynamical zeta functions  which are, unlike the Selberg- or the
Artin-Weil zeta functions, not related to arithmetics.

In the present paper we address this problem for the Ruelle zeta function of a certain
subshift of finite  type which in the physics literature has become known as the Kac-Baker model
(see \cite{Ka66}).  M.~Kac got interested in this model while trying to understand the
mathematics behind the phenomenon of phase transitions in systems
of statistical mechanics with weak long-range interactions like in the van
der Waals gas (see \cite{Ba61}).  The model he considered is an Ising spin system on a 1-dimensional
lattice with a 2-body interaction given by a finite superposition of terms
decaying exponentially fast with the distance between the different spins.
His real interest was certainly in a system with a  continuous superposition of such
exponentially decaying terms to model also interactions decaying only polynomially fast as
is the case in the van der Waals gas. However his method did not allow him to treat this
limiting case in a rigorous way.

We have chosen the Kac model since its dynamical zeta function can be
understood rather well by the transfer operator method. On the other hand there seems to
be no obvious connection of this zeta function to any arithmetic zeta functions for which
a general Riemann hypothesis is known to hold.

There exist two rather different transfer operators for this model
which allow to express its zeta function as Fredholm determinants of these operators 
(see \cite{Ma80}, \cite{ViMa77}, \cite{Ka66}).
Up to now, however, it was not known how these two operators are related to each other.
Our investigations show that the Ruelle operator is basically equivalent
to Kac's original transfer operator through a Segal-Bargmann transformation
establishing a unitary map between the Hilbert space of square integrable
functions on the real line, where the Kac operator acts, and the Fock space of entire
functions on the complex plane square integrable with respect
to a certain weight function, to which the Ruelle operator can be restricted.
Since the Kac operator $\mathcal{K}_{\beta }$ is a symmetric, positive definite trace class operator for positive $\beta$ and has real spectrum also for negative $\beta$, and the Ruelle operator $\mathcal{L}_{\beta}$ defines a  family of trace class operators holomorphic in the variable $\beta$  we can show that the zeta functions of a whole
class of Kac models extend meromorphically to the entire complex
$\beta$-plane and have infinitely many nontrivial zeros on the real line.
For a special case of the parameters we can also show the existence of infinitely many
``trivial''  zeros of this dynamical zeta function located on a line parallel to the
imaginary axis in the complex $\beta$-plane. Thus for this function an analogue
of the  Riemann hypothesis seems plausible.
The present paper generalizes analogous results in \cite{HiMa01}
for the case of an interaction consisting of a single exponentially decaying term.

In detail the present paper is organized as follows: in Section
\ref{ROKM} we recall the definition of the Kac-Baker models and
derive their Ruelle transfer operators. Further, we show how the
dynamical zeta function of these models can be expressed through
Fredholm determinants of the Ruelle operators. In Section
\ref{KGO} we derive,  basically following Gutzwiller (see
\cite{Gu82}), the Kac operator appropriate to our problem and show
that the zeta function can be expressed for positive $\beta$ also
in terms of Fredholm determinants of this Kac operator. In Section
\ref{HFMF} we show how the kernel of the Kac operator is related
to a certain form of Mehler's formula for the Hermite functions
which allows us to diagonalize an integral operator closely
related to the Kac operator. In Section \ref{FSSBT} we introduce
the Fock spaces and the Segal-Bargmann transform and show how the
Ruelle and Kac operators can be directly related to each other.
There we show that for real $\beta$ the two operators have the
same spectrum and give explicit expressions relating the
eigenfunctions of the two operators for nonvanishing eigenvalues.
In Section \ref{ZRRF} we derive the asymptotic behavior of the
eigenvalues of the  Ruelle operator for large positive and
negative values of $\beta$ and apply it to get the results on the
location of poles and zeros of the zeta function on the real line.
In Section \ref{MatEl} we give explicit expressions for the matrix
elements of a modified Kac-Gutzwiller operator in the Hilbert
space basis given by the Hermite functions which seem best suited
for future numerical calculations.

\section{The Ruelle operator for the Kac-Baker model}
\label{ROKM}

The generalized Kac model describes a 1-dimensional
lattice spin system with a 2-body interaction which is a superposition
of finitely many exponentially decaying terms.
More precisely, for $F:=\{\pm 1\}$,
$\underline{\xi}=(\xi_n)_{n\in \Z_+}\in F^{\Z_+}$, and
 $i,j\in \Z_+=\left\{0,1,2,\ldots\right\}$ we set
\[
\phi_{ij}(\underline{\xi})
:=-\xi_{i}\xi_{j}\sum _{l=1}^{m}J_{l}\lambda ^{\left| i-j\right| }_{l},
\]
where $m\in \N$ and the parameters  $J_{l}>0$ and $0<\lambda
_{l}<1$ are fixed and describe the interaction strengths and the
different decay rates. The interaction energy of a configuration
$\underline{\xi} \in F^{\mathbb {Z}_{+}}$ when restricted to the
finite sublattice $\mathbb {Z}_{\left[ n-1\right] }=\left\{
0,1,\ldots ,n-1\right\}$ is then given as
\[
U_{n}(\underline{\xi} )
:=U_{\mathbb {Z}_{\left[n-1\right]}}(\underline{\xi} )
=\sum ^{n-1}_{i=0}\sum ^{\infty }_{j=1}\phi_{i,i+j}(\underline{\xi}).
\]
When inserting the explicit form of $\phi$ one gets
\[
U_{n}(\underline{\xi})
=-\sum ^{n-1}_{i=0}\sum ^{\infty }_{j=1}\xi _{i}\xi _{i+j}
  \sum ^{m}_{l=1}J_{l}\lambda ^{j}_{l}.
\]
For $\beta\in \C$ the partition functions $Z_{n}(\beta )$  for the finite sublattices
$\Z_{\left[n-1\right]}$
with periodic boundary conditions are defined as
\begin{equation}\label{Zndef}
Z_{n}(\beta )
:=\sum _{\underline{\xi} \in {\Per}_{n}}\exp\left( -\beta U_{n}(\underline{\xi} )\right),
\end{equation}
where ${\Per}_{n}$ denotes the set of configurations
$\underline{\xi }\in F^{\Z_{+}}$
which are periodic with period $n$. That means
${\underline{\xi} }\in {\Per}_{n}$ if and only if $\xi _{i+n}=\xi _{i}$ for all
$i\in \Z_{+}$.
Defining  the shift   $\tau :F^{\Z_{+}}\rightarrow F^{\Z_{+}}$
by
\[
\left( \tau \underline{\xi }\right) _{i}:=\xi _{i+1}\:
\quad  \mbox{if } {\underline{\xi}}=\left( \xi _{i}\right) _{i\in \Z_+},
\]
one has
$\Per_{n}=\Fix\,\tau ^{n}
=\left\{ \underline{\xi}\in F^{\Z_{+}}:\tau ^{n}{\underline{\xi}}
={\underline{\xi}}\right\}
$.

To the dynamical system
 \( \left( F^{\Z_{+}},\tau \right)  \)
one can associate the Ruelle zeta function
\begin{equation}
\label{RuelleZetadef}
\zeta_{R}(z,\beta ):=\exp\left( \sum ^{\infty }_{n=1}\frac{z^{n}}{n}Z_{n}(\beta )\right).
\end{equation}
Note that $|U_n(\underline{\xi})|\le
n\sum_{l=1}^m\frac{\lambda_lJ_l}{1-\lambda_l}=:nc$ so that
$|Z_n(\beta)|\le (2e^{|\beta|c})^n$. Therefore the series defining
$\zeta_R$ converges in a neighborhood of $(0,0)$ in $\C^2$. We
will show in Proposition \ref{zetaRmero} that $\zeta_R$ can in
fact be extended to a meromorphic function on $\C^2$.

To determine the analytic properties of this function one makes use of
the transfer operator technique. Note first that the configuration space
$F^{\Z_+}$ is  compact and metrizable with respect to the product topology.
For each $\beta\in \C$ one can define the Ruelle transfer operator $\mathcal{L}_\beta$
which will act on the space of observables of the lattice spin system,
i.e. on $\mathcal{C}(F^{\Z_{+}})$, the space
of continuous  functions on  $F^{\mathbb {Z_{+}}}$.
In this framework one sets
\[
\left( \mathcal{L}_\beta f\right) \left( {\underline{\xi}}\right)
:=\sum _{{\underline{\eta}}\in \tau ^{-1}(\underline{\xi})}
   \exp\left( -\beta U_{1}( {\underline{\eta}})\right)f({\underline{\eta}}).
\]
Inserting the explicit expression for $U_{1}(\underline{\eta})$
one finds
\begin{equation}
\label{RuelleTransfer}
\left( {\mathcal{L}}_{\beta }f\right) (\underline{\xi})
=\sum _{\sigma =\pm 1}\exp\left( \beta \sigma  \sum ^{\infty }_{j=0}\xi _{j}
 \sum ^{m}_{l=1}J_{l}\lambda ^{j+1}_{l}\right)
  f\left( \left( \sigma ,\underline{\xi}\right)\right),
\end{equation}
where $(\sigma,\underline{\xi}):=\underline{\eta}$ with $\eta_0=\sigma$ and
$\eta_j=\xi_{j-1}$ for all $j\in \N$.
Generalizing the arguments for the case $m=1$ in \cite{Ma80}  one introduces the map
$\underline{z}=\left(z_{l}\right)_{l=1,..,m}\colon F^{\Z_+}\to\R^m$ defined by
\[
z_{l}(\underline{\xi})
:=\sum ^{\infty }_{i=0}\xi _{i}\lambda ^{i+1}_{l}.
\]
Since
$z_{l}(  \sigma ,\underline{\xi})
=\sigma \lambda _{l}+\lambda _{l}z_{l}( {\underline{\xi}})
$
for all $l$, the operator $\mathcal{L}_{\beta }$ leaves the space of functions
$\underline{z}^*(\varphi):=\varphi\circ \underline{z}$ with $\varphi\in \mathcal{C}(\R^m)$
invariant. Thus we obtain a factorization
$$
\xymatrix{
\mathcal{C}(F^{\Z_+})\ar[r]^{\cL_\beta}&\mathcal{C}(F^{\Z_+})\\
\mathcal{C}(\R^m)\ar@{.>}[r]\ar[u]^{\underline{z}^*}&\mathcal{C}(\R^m)\ar[u]_{\underline{z}^*}
}
$$
of $\cL_\beta$ through $\mathcal{C}(\R^m)$, which we will denote
again by $\mathcal{L}_{\beta }$ and which  is of the form

 \begin{equation}\label{RuelleTO} \left(
\mathcal{L}_{\beta }g\right) ( \underline{z})
=e^{\beta\underline{J}\cdot\underline{z}}g(\boldsymbol{\Lambda}\underline{z}
   +\underline{\lambda})
+e^{-\beta\underline{J}\cdot\underline{z}}g(\boldsymbol{\Lambda}\underline{z}
   -\underline{\lambda}),
\end{equation} where $\boldsymbol{\Lambda}$ is the diagonal matrix
with diagonal elements $\lambda_1,\ldots,\lambda_m$ and
$\underline{J}=(J_1,\ldots,J_m)$. 

 Indeed, by the
Ruelle-Perron-Frobenius Theorem (see \cite[Chap.4]{Zi00}) the
iterates under $\mathcal{L}_{\beta}$ of the constant function
$f(\underline{\xi})=1$ in $\mathcal{C}(F^{\Z_+})$ converge
uniformly to the eigenfunction belonging to the leading eigenvalue
of $\mathcal{L}_{\beta}$ and hence this eigenfunction belongs to
the space $\mathcal{C}(\R^{m})$. Therefore, when restricting the
operator to the space $\mathcal{C}(\R^m)$ one does not lose the
leading eigenvalue which is the most important one from the
physical point of view. The physically most satisfying operator is
obtained in fact by restricting the domain of the operator still
further.

>From the form of the operator $\mathcal{L}_{\beta}$ in the space
$\mathcal{C}(\R^{m})$ we see that, if $g$  is a holomorphic function on the polycylinder
\[
D=\left\{\underline{z}\in \C^m:\left| z_{l}\right| <R_{l},l=1,\ldots,m\right\} \]
with $R_{l}>\frac{\lambda _{l}}{1-\lambda _{l}}$, then also
$\left( {\mathcal{L}}_{\beta }g\right) \left(\underline{z}\right)$ is such a function.
In fact, much more is true: Denote by
$\mathcal{B}(D)$ the Banach space of  holomorphic functions on
$D$ which extend continuously to the closure $\oline{D}$ of $D$.
Then, according to \cite[Appendix B]{Ma80},  we have

\begin{Prop}\label{RTOnuclear}
${\mathcal{L}}_{\beta }\colon \mathcal{B}(D)\to \mathcal{B}(D)$
is with respect to the parameter $\beta\in\C$ a holomorphic family of  nuclear operators of
degree $0$ in the sense of Grothendieck
(see \cite{Gr55}). In particular all the $\cL_\beta$ are of trace class.
\qed
\end{Prop}
The trace can be computed using the holomorphic version
of the Atiyah-Bott fixed point formula (see \cite{AtBo67} and \cite{Ma80}, Appendix B
or \cite{Ru94},\S 1.12):

\begin{Lemma}\label{traceformula}
Fix $\varphi\in \mathcal{B}(D)$ and a continuous map $\psi:\oline{D}\to D$ which is
holomorphic on $D$.
Then $\psi$ has a unique fixed point $\underline{z}^\fix\in D$ and the composition operator
\(A: \mathcal{B}(D) \to \mathcal{B}(D)\) defined by  $Ag:= \varphi\cdot (g\circ\psi)$
is trace class with trace
\[
\trace(A)= \frac{\varphi(\underline{z}^\fix)}{\det\big(1-\psi'(\underline{z}^\fix)\big)}.
\]
\qed
\end{Lemma}

\begin{Prop}\label{Zntrace}
The partition function $Z_{n}(\beta )$ of the Kac-Baker model can be expressed through the
traces of the powers of the Ruelle transfer operator ${\mathcal{L}}_{\beta }$
via
\[
Z_{n}(\beta )=\left(\prod ^{m}_{l=1}\left( 1-\lambda ^{n}_{l}\right)\right)\trace \cL_\beta^n.
\]
\end{Prop}

\pf
The defining equation (\ref{Zndef}) for the partition function $Z_{n}(\beta)$
can be rewritten as
\[
Z_{n}(\beta)
=\sum_{\underline{\xi}\in \Per_{n}}
     \exp\left(\beta\sum_{l=1}^{m}J_{l}\sum_{k=0}^{n-1}\sum_{i=1}^{\infty}
                \xi_{k}\xi_{k+i}\lambda_{l}^{i}\right).
\]
Using the fact that  $\underline{\xi}\in\Per_n$ implies
$\xi_{i+n}=\xi_{i}$  for all $i\in\mathbb{Z}_{+}$ one gets
\[
Z_{n}(\beta)
= \sum_{\underline{\sigma} \in F^{n}}
    \exp\left( \beta \sum_{l=1} ^{m} \frac {J_{l}} {1-\lambda_{l}^{n}}
              \sum_{k=0}^{n-1} \sum_{i=1}^{n} \sigma_{k} \sigma_{k+i}\lambda_{l}^{i}\right),
\] where \(\sigma_{n+i}=\sigma_{i}\) for all \(i\).
On the other hand the $n$-th iterate of the transfer operator
$\mathcal{L}_{\beta}$ from (\ref{RuelleTO}) acting on the Banach space $\mathcal{B}(D)$ is given by
\begin{eqnarray*}
&&\hskip -2em
\left(\mathcal{L}_{\beta}^{n}g\right)(\underline{z}) =\\
&=&\sum_{\underline{\sigma} \in F^{n}} \exp \left( \beta\sum_{l=1}^{m}J_{l}
    \left( \sum_{k=1}^{n} \sigma_{k} \lambda_{l}^{n-k} z_{l} + \sum_{k=1}^{n-1}
           \sum_{i=1}^{n-k} \sigma_{k} \sigma_{k+i} \lambda_{l}^{i}
    \right)\right)
    g\left(\sum_{i=1}^{n}\sigma_{i} \underline{\lambda}^{i}
    + \boldsymbol{\Lambda}^{n} \underline{z}\right),
\end{eqnarray*}
where $\underline{\lambda}^i:=(\lambda_1^i,\ldots,\lambda_m^i)$.
 We apply Lemma \ref{traceformula} to the maps
$\psi_{\underline{\sigma}}$ defined by
\[
\psi_{\underline{\sigma}}(\underline{z}) := \left( \sum_{i=1}^{n}
\sigma_{i} \underline{\lambda}^{i} + \boldsymbol{\Lambda}^{n}
\underline{z}\right)
\]
for which the fixed points are given by
\[
\underline{z}_{\underline{\sigma}}^\fix =
(1-\boldsymbol{\Lambda}^n)^{-1}\sum_{i=1}^{n} \sigma_{i}
\underline{\lambda}^{i}.
\]

 The result is
\begin{eqnarray*}
\trace{\mathcal{L}_{\beta}^{n}}
&=& \frac{1}{\prod _{l=1}^{m}(1-\lambda_{l}^{n})} \sum_{\underline{\sigma}\in F^{n}}
    \exp\left( \beta \sum_{l=1}^{m}\frac {J_{l}} {1-\lambda_{l}^{n}}
    \left( \sum_{k=1}^{n} \sum_{i=1}^{n} \sigma_{k} \sigma_{i}
    \lambda_{l}^{n-k+i} +\right.\right.\\
&&  \phantom{\frac{1}{\prod _{l=1}^{m}(1-\lambda_{l}^{n})} \sum_{\underline{\sigma}\in F^{n}}
    \exp( } \left.\left.+    \sum_{k=1}^{n-1} \sum_{i=1}^{n-k} \sigma_{k} \sigma_{k+i} \lambda_{l}^{i} -
    \sum_{k=1}^{n-1} \sum_{i=1}^{n-k} \sigma_{k} \sigma_{k+i} \lambda_{l}^{n+i}
     \right)\right).
\end{eqnarray*}
But
\[
\sum_{k=1}^{n-1} \sum_{i=k+1}^{n} \sigma_{k} \sigma_{i} \lambda_{l}^{n-k+i}
= \sum_{k=1}^{n-1} \sum_{i=1}^{n-k} \sigma_{k} \sigma_{k+i} \lambda_{l}^{n+i}
\]
and hence
\begin{eqnarray*}
\trace \mathcal{L}_{\beta}^{n}
&=& \frac{1}{\prod_{l=1}^{m} (1-\lambda_{l}^{n})}
  \sum_{\underline{\sigma}\in F^{n}}
  \exp\left( \beta\sum_{l=1}^{m}
             \frac{J_{l}}{1-\lambda_{l}^{n}}
             \left( \sigma_{n} \sum_{i=1}^{n} \sigma_{i} \lambda_{l}^ {i} +
                   \sum_{k=1}^{n-1} \sum_{i=1}^{n-k} \sigma_{k} \sigma_{k+i} \lambda_{l}^{i}
+\right.\right.\\
&&\phantom{\frac{1}{\prod_{l=1}^{m} (1-\lambda_{l}^{n})}
  \sum_{\underline{\sigma}\in F^{n}}
  \exp\left( \beta\sum_{l=1}^{m}
             \frac{J_{l}}{1-\lambda_{l}^{n}}\right.}
             \left.\left.
                \,\,+ \sum_{k=1}^{n-1} \sum_{i=1}^{k} \sigma_{k} \sigma_{i} \lambda_{l}^{n-k+i}
             \right)
     \right).
\end{eqnarray*}
Changing the order of summation in the last double sum we finally get
\[
\trace{\mathcal{L}_{\beta}^{n}}
= \frac{1}{\prod_{l=1}^{m}(1-\lambda_{l}^{n})} \sum_{\underline{\sigma}\in F^{n}}
  \exp\left(\beta\sum_{l=1}^{m} \frac{J_{l}}{1-\lambda_{l}^{n}}
            \sum_{k=1}^{n} \sum_{i=1}^{n} \sigma_{k} \sigma_{k+i}\lambda_{l}^{i}
       \right)
\]
which up to the factor $\prod_{l=1}^{m}\frac{1}{1-\lambda_{l}^{n}}$ is just the partition
function $Z_{n}(\beta)$.
\qed

In view of the identity
$\prod ^{m}_{l=1}\left( 1-\lambda ^{n}_{l}\right)
 =\sum _{\underline{\alpha} \in \left\{ 0,1\right\} ^{m}}
  \left( -1\right) ^{\left|\underline {\alpha }\right| }
  \lambda ^{n\alpha _{l}}_{1}\cdots \lambda ^{n\alpha _{m}}_{m}$
Proposition \ref{Zntrace} yields
\[
Z_{n}(\beta )
=\sum _{\underline{\alpha }\in \left\{ 0,1\right\} ^{m}}
       \left(-1\right)^{\left|\underline{\alpha}\right|}
              \trace\left( \left( \prod ^{m}_{l=1}\lambda ^{\alpha _{l}}_{l}\right) ^{n}
             \mathcal{L}^{n}_{\beta }\right),
\]
so that
\begin{eqnarray*}
\zeta _{R}\left( z,\beta \right)
&=&\exp \left(\sum ^{\infty }_{n=1}\frac{z^{n}}{n}
   \sum _{\underline{\alpha }\in \left\{ 0,1\right\} ^{m}}\left(-1\right)^{\left|\underline{\alpha}\right|}
   \trace\left( \left( \prod ^{m}_{l=1}\lambda ^{\alpha _{l}}_{l}\right)
   \cL_{\beta }\right) ^{n}\right) \\
&=&\exp\left(\trace \sum _{\underline{\alpha }\in \left\{ 0,1\right\} ^{m}}
   (-1)^{\left|\underline{\alpha }\right| }(-1)
   \log \left( 1-z\left( \prod ^{m}_{l=1}\lambda ^{\alpha _{l}}_{l}\right)
    \cL_{\beta }\right)\right) \\
&=&\prod _{\underline{\alpha }\in \left\{ 0,1\right\} ^{m}}
   \det \left( 1-z\prod ^{m}_{l=1}\left( \lambda ^{\alpha _{l}}_{l}\right)
   \cL_{\beta }\right) ^{\left( -1\right) ^{\left|\underline{\alpha }\right| +1}}.
\end{eqnarray*}

Together with Proposition \ref{RTOnuclear}
this proves the following  proposition.

\begin{Prop}\label{zetaRmero}
The Ruelle zeta function
$\zeta_{R}(z,\beta )$ defined in (\ref{RuelleZetadef})
for the Kac-Baker model with decay rates
$\underline{\lambda}=(\lambda_1,\ldots,\lambda_m)\in]0,1[^m$
can be extended to a meromorphic family
$\left(\zeta_R(\cdot ,\beta)\right)_{\beta\in\C}$ of meromorphic functions
on $\C$ via the formula
$$\zeta _{R}\left( z,\beta \right)=
\prod _{\underline{\alpha }\in \left\{ 0,1\right\} ^{m}}
   \det \left( 1-z \underline{\lambda}^{\underline{\alpha}}
   \cL_{\beta }\right) ^{\left( -1\right) ^{\left|\underline{\alpha }\right| +1}}.
$$
\qed
\end{Prop}

\begin{Rem}
The analytic properties in the variable \(z\) are well known for general lattice spin systems with exponentially
fast decaying interactions.
For $m=1$ Proposition \ref{zetaRmero} was proved in \cite{Ma80}.
\qed
\end{Rem}

\section{The Kac-Gutzwiller operator}
\label{KGO}

In  \cite{Ka59}  M.\,Kac found another operator whose traces are directly
related to the partition functions
 \( Z_{n}(\beta ) \)
 of the Kac-Baker model.
 Kac did not work with periodic but open boundary conditions and hence
 his operator has to be modified a bit to give the partition functions we use here.
 In the case
\( m=1 \)
 M.~Gutzwiller  already derived this operator (see \cite{Gu82})
and the case of general
 \( m\geqslant 1 \)
 can be handled similarly.
 Just like Gutzwiller and Kac we start with an identity for Gaussian integrals known
  already to C.\,Cram\'er: For a positive definite $(n\times n)$-matrix
 \(\mathbb{A} \) and $\underline{x}\in \R^n$  we have (see \cite{Cr46} or the
 appendix of \cite{Fo89})
\begin{equation}\label{Cramerformel}
e^{\frac{1}{2}\left(\underline{x}\cdot\mathbb{A}\underline{x}\right)}
=(2\pi )^{-\frac{n}{2}}(\det\mathbb{B})^{\frac{1}{2}}\int _{\R^n}
e^{\underline{x}\cdot\underline{z}} e^{-\frac{1}{2}
 (\underline{z}\cdot\mathbb{B}\underline{z})} d\underline{z},
\end{equation}
 where \(\mathbb{B}=\mathbb{A}^{-1} \).

 Any periodic configuration $\underline{\xi }^{p}\in \Per_{n}$
can be extended to a periodic configuration on the entire lattice
$\Z$ with the same period which we again denote  by
$\underline{\xi}^{p}$.

Now consider the $(n\times n)$-matrix $\mathbb{A}^{(l)}$ given by
$$ \mathbb {A}^{(l)}_{i,j}
 =  \beta J_{l}\sum ^{+\infty }_{k=-\infty }\exp\left(-\gamma_{l} \left| i-j+nk\right|\right),
    \quad 0\leq i,j\leq n-1,$$
where we choose the constants $\gamma_l$  so that
$e^{-\gamma_l}=\lambda_l$ with the coupling constants $\lambda_l$
of the Kac-Baker model. Using $\underline{\xi }_{n}^{p}  =  \left(
\xi^p _{0},\ldots ,\xi^p _{n-1}\right) $ one finds (see \cite[\S
6]{Gu82} for this and the following results on matrix
calculations)
\[
\frac{\beta }{2}J_{l}\sum ^{n-1}_{i=0}\sum ^{+\infty }_{j=-\infty }\xi ^{p}_{i}
     \xi ^{p}_{j}\exp\left( -\gamma _{l}\left| i-j\right|\right)
=\frac{1}{2}\left(\underline{\xi }_{n}^{p}\cdot \mathbb{A}^{(l)}
  \underline{\xi }_{n}^{p}\right).
\] But
\begin{equation}\label{IDT}
\frac{\beta }{2}\sum\limits ^{m}_{l=1}J_{l}
  \sum\limits ^{n-1}_{i=0}\sum\limits ^{+\infty }_{j=-\infty }
  \xi ^{p}_{i}\xi ^{p}_{j}\lambda ^{\left| i-j\right| }_{l}
= -\beta U_{n}\left(\underline {\xi }^p\right) +\frac{\beta }{2}n
\sum\limits ^{m}_{l=1}J_{l},
\end{equation}
i.e. the left hand side is, up to
the constant term given by the sum of the $J_{l}$,
 just the interaction energy of the periodic configuration
 \(\underline{\xi }^p\in \Per_{n} \)
 for the Kac model.
\newline
The matrix
 \(\mathbb {B}^{(l)}=\left(\mathbb{A}^{(l)}\right)^{-1} \)
 has the following form
\[
\mathbb {B}^{(l)}=\frac{1}{\beta J_{l}\sinh \gamma _{l}}\,\begin{pmatrix}
\cosh \gamma _{l}& -\frac{1}{2}    & 0               &\ldots      &\ldots           & 0               & -\frac{1}{2}\\
-\frac{1}{2}     &\cosh \gamma _{l}& -\frac{1}{2}    & 0          &\ldots           & 0               &  0          \\
0                &-\frac{1}{2}     &\cosh \gamma _{l}&-\frac{1}{2}&0                &\ldots           &  0          \\
\vdots           & \ddots          & \ddots          &\ddots      &\ddots           &\ddots           &\vdots       \\
0                &\ldots           &0                &-\frac{1}{2}&\cosh \gamma _{l}&-\frac{1}{2}     &0            \\
0                &0                &\ldots           &0           &-\frac{1}{2}     &\cosh \gamma _{l}&-\frac{1}{2} \\
-\frac{1}{2}     &0                &\ldots           &\ldots      & 0               & -\frac{1}{2}    &\cosh \gamma _{l}
\end{pmatrix},
\]
which for positive $\beta$ is a positive definite matrix
with determinant
\[
\det \mathbb {B}^{(l)}=\frac{4}{\left( 2\beta J_{l}\sinh \gamma _{l}\right) ^{n}}
\left( \sinh \frac{n\gamma }{2}\right) ^{2}.
\]
That $\mathbb{B}^{(l)}$ is indeed positive definite one can see as follows:
for arbitrary $\underline{x}\in \mathbb{R}^{n}$ one finds
\[\left(\underline{x}\cdot\mathbb{B}^{(l)}\underline{x}\right)=\frac{1}{\beta J_{l}}\left(\sum_{i=1}^{n}\coth(\gamma_{l}) x_{i}^{2}-\frac{\sum^{n}_{i=1}x_{i}x_{i-1}}{\sinh \gamma_{l}}\right)\]where $x_{0}= x_{n}.$
A simple calculation (see also Proposition \ref{KacKernProp})
shows that the following identity holds
\begin{equation}\label{FORM}
  \coth (\gamma _{l})\sum ^{n}_{i=1}x^{(l)2}_{i}
 -\frac{\sum\limits ^{n}_{i=1}x^{(l)}_{i}x^{(l)}_{i-1}}{\sinh \gamma _{l}} = \frac{1}{2}\left(\tanh \left(\frac{\gamma _{l}}{2}\right)\sum ^{n}_{i=1}\left( x^{(l)2}_{i}+x^{(l)2}_{i-1}
     \right)
+\frac{\sum ^{n}_{i=1}\left( x^{(l)}_{i}-x^{(l)}_{i-1}\right) ^{2}}{\sinh \gamma _{l}}\right)
\end {equation}
and hence for $\beta J_{l}>0$ and $\gamma_{l}>0$ one finds \(\left(\underline{x},\mathbb{B}^{(l)}\underline{x}\right) \geqslant 0\).
Inserting (\ref{IDT}) into formula (\ref{Cramerformel}) one calculates

\begin{eqnarray*}
e^{ -\beta U_{n}\left( \underline{\xi }^{p}\right)
            +\frac{n\beta }{2}\sum\limits ^{m}_{l=1}J_{l}}
&=&\prod\limits ^{m}_{l=1} e^{\frac{1}{2}
   \left(\underline{\xi }_{n}^{p}\cdot \mathbb{A}^{(l)}\underline{\xi }_{n}^{p}\right)}\\
&=&\prod\limits ^{m}_{l=1}\left(
   ( 2\pi )^{-\frac{n}{2}}\frac{2}{\left( 2\beta J_{l}\sinh \gamma _{l}\right) ^{\frac{n}{2}}}
   \sinh \left(\frac{n\gamma _{l}}{2}\right)\int\limits _{\R^n}
   e^{ -\frac{1}{2}
   \left(\underline{z}^{\left( l\right) }\cdot \mathbb{B}^{(l)}\underline{z}^{(l)}\right)}
   e^{\underline{\xi }_n^{p}\cdot\underline{z}^{\left(l\right)}}d\underline{z}^{\left(l\right)}\right)\\
&=&\frac{2^{m}}{2^{nm}}\left( \prod\limits ^{m}_{l=1}
   \frac{\sinh \frac{n\gamma _{l}}{2}}{\left( \beta \pi J_{l}
    \sinh \gamma _{l}\right)^{\frac{n}{2}}}\right)
    \int\limits_{\R^{n}}\ldots \int\limits _{\R^{n}}
    e^{ -\frac{1}{2}\sum\limits ^{m}_{l=1}\left(\underline {z}^{(l)}\cdot\mathbb {B}^{(l)}\underline{z}^{(l)}\right)}
    e^{\sum\limits ^{m}_{l=1}\left( \underline{\xi }_n^{p}\cdot\underline{z}^{(l)}\right)}
    d\underline{z}^{(1)}\ldots d\underline{z}^{(m)}.
\end{eqnarray*}
The change $\underline
{z}^{(l)}=:\sqrt{\beta J_{l}}\, \underline{x}^{(l)}$
of  integration variables yields
\begin{eqnarray*}
&&\hskip -2em
e^{-\beta U_{n}(\underline{\xi}^{p} )+\frac{n\beta }{2}\sum ^{m}_{l=1}J_{l}}=\\
&=&\frac{2^{m}}{2^{nm}}\prod ^{m}_{l=1}
    \frac{\left(\sinh \frac{n\gamma _{l}}{2}\right)\left( \beta J_{l}\right) ^{\frac{n}{2}}}
    {\left( \beta J_{l}\pi \sinh \gamma _{l}\right) ^{\frac{n}{2}}}
    \int\limits _{\R^{n}}\ldots \int\limits _{\R^{n}}
  e^{-\frac{1}{2}\sum ^{m}_{l=1}\beta J_{l}\left(\underline{x}^{(l)}\cdot \mathbb{B}^{(l)}\underline{x}^{(l)}\right)}
  e^{\sum ^{m}_{l=1}\sqrt{\beta J_{l}}\left(\underline{\xi }_n^{p}\cdot\underline{x}^{(l)}\right)}
 d\underline{x}^{(1)}\ldots d\underline{x}^{(m)}.
\end{eqnarray*}
Performing the summation over all the periodic configurations
$\underline {\xi }\in \Per_{n}$  then gives
\begin{eqnarray*}
e^{\frac{n\beta}{2} \sum\limits ^{m}_{l=1}J_{l}}
   \sum _{\underline{\xi }\in \Per_{n}}e^{-\beta U_{n}
   \left(\underline{\xi }\right)}
&=&\frac{2^{m}}{2^{n(m-1)}}\prod ^{m}_{l=1}\left(
    \frac{\sinh \frac{n\gamma _{l}}{2}}{(\pi \sinh \gamma _{l})^{\frac{n}{2}}}
    \right)\cdot\\
&&\cdot    \int\limits _{\R^{n}}\ldots \int\limits _{\R^{n}}
    e^{-\frac{1}{2}\sum ^{m}_{l=1}\beta J_{l}\left(\underline{x}^{(l)}\cdot\mathbb{B}^{(l)}\underline{x}^{(l)}\right)}
    \prod ^{n}_{i=1}\cosh \left( \sum ^{m}_{l=1}\sqrt{\beta J_{l}}x^{(l)}_{i}\right)
    d\underline{x}^{(1)}\ldots d\underline{x}^{(m)}.
\end{eqnarray*}
 Inserting the explicit form of the matrix
\( \mathbb{B}^{(l)}=\left(\mathbb{A}^{(l)}\right)^{-1} \)
 we find
\begin{eqnarray*}
e^{\frac{n\beta }{2}\sum ^{m}_{l=1}J_{l}}\, Z_{n}(\beta )
&=&\frac{2^{m}}{2^{n(m-1)}}\prod ^{m}_{l=1}
   \frac{\sinh \frac{n\gamma _{l}}{2}}{\left( \pi \sinh \gamma _{l}\right) ^{\frac{n}{2}}}
   \cdot\\
& & \cdot
   \int\limits _{\R^{n}}\ldots \int\limits _{\R^{n}}
   \exp\left(-\frac{1}{2}\sum ^{m}_{l=1}
       \left( (\coth \gamma _{l})\sum ^{n}_{i=1}(x_{i}^{(l)})^{2}
   -\frac{\sum\limits ^{n}_{i=1}x^{(l)}_{i}x^{(l)}_{i-1}}{\sinh \gamma _{l}}\right)\right)\cdot\\
& &\cdot
   \prod ^{n}_{i=1}\cosh \left( \sum ^{m}_{l=1}\sqrt{\beta J_{l}}x^{(l)}_{i}\right)
d\underline{x}^{(1)}\ldots d\underline{x}^{(m)},
\end{eqnarray*}
 where
 \( x_{0}^{(l)}=x_{n}^{(l)}\).

For $\beta \geqslant 0$ we introduce the kernel function
\begin{equation}\label{KacGutzKern}
\mathcal{K}_{\beta }\left(\underline{\xi},\underline{\eta}\right)
:=\textstyle{
\frac{\left( \cosh \left( \sum\limits ^{m}_{l=1}\sqrt{\beta J_{l}}\xi_{l}\right)
  \cosh \left( \sum\limits ^{m}_{l=1}\sqrt{\beta J_{l}}\eta_{l}\right) \right) ^{\frac{1}{2}}}
  {2^{(m-1)}\prod\limits ^{m}_{l=1}\left( \pi \sinh \gamma _{l}\right) ^{\frac{1}{2}}}
  \exp\left(-\frac{1}{4}\left(\sum ^{m}_{l=1}\left(( \tanh \frac{\gamma _{l}}{2})
  \left( \xi^{2}_{l}+\eta^{2}_{l}\right) +\frac{\left( \xi_{l}-\eta_{l}\right) ^{2}}
  {\sinh \gamma _{l}}\right)\right)\right),
}
\end{equation}
where
$\underline{\xi}=\left( \xi_{1},\ldots ,\xi_{m}\right) \in \R^{m}$ and
$\underline{\eta}=\left( \eta_{1},\ldots ,\eta_{m}\right) \in \R^{m}$.
We call the  associated operator $\mathcal{K}_\beta$ on $L^2(\R^{m},d\underline{\xi})$ defined by
\[
\left(\mathcal{K}_\beta f\right) \left( \underline{\xi}\right)
=\int\limits _{\R^{m}}\mathcal{K}_\beta\left(\underline{\xi},\underline{\eta}\right)
 f\left( \underline{\eta}\right) d\underline{\eta}
\]
the Kac-Gutzwiller transfer operator. Note that the kernel of $\cK_\beta$ decreases fast enough
to ensure that $\cK_\beta$ is of trace class with
$\trace\cK_\beta=\int_{\R^m}\cK_\beta(\underline{\xi},\underline{\xi})d\underline{\xi}$.
Calculating the trace of the iterates of $\cK_\beta$ and comparing the result to
the above formula for the partition function \(Z_{n}(\beta)\) after inserting (\ref{FORM}) we find
\begin{equation}\label{ZnformelK}
Z_{n}(\beta )=2^{m}\prod ^{m}_{l=1}\left( \sinh \frac{n\gamma _{l}}{2}\right)
e^{-\frac{n\beta }{2}\sum ^{m}_{l=1}J_{l}}\trace \mathcal{K}^{n}_{\beta }.
\end{equation}

To simplify the expression for
$Z_{n}(\beta )$ we  introduce the rescaled Kac-Gutzwiller operator
$\mathcal{G}_{\beta }:L^{2}\left( \R^{m},d\underline{\xi}\right) \to L^{2}
\left( \R^{m},d\underline{\xi}\right)$
defined by
\begin{equation}\label{KacGutzmod}
\mathcal{G}_{\beta }:=\prod ^{m}_{l=1}\left( \lambda _{l}e^{\beta J_{l}}\right) ^{-\frac{1}{2}}
\mathcal{K}_{\beta }.
\end{equation}
In view of
\[ Z_{n}(\beta )=(\prod\limits ^{m}_{l=1}\left( 1-\lambda ^{n}_{l}\right)) \trace
(\prod\limits ^{m}_{l=1}e^{\frac{n\gamma _{l}}{2}}e^{-\frac{\beta J_{l}n}{2}})\mathcal{K}^{n}_{\beta }, \]
which is a simple reformulation of (\ref{ZnformelK})
we finally have shown the following proposition.

\begin{Prop}\label{Zntrace2}
For $\beta\geqslant 0$ the partition function $Z_{n}(\beta )$ of the Kac-Baker model can be expressed through the
traces of the powers of the rescaled Kac-Gutzwiller operator ${\mathcal{G}}_{\beta }$
via
\[
Z_{n}(\beta )=\left(\prod ^{m}_{l=1}\left( 1-\lambda ^{n}_{l}\right)\right)\trace \cG_\beta^n.
\]
\qed
\end{Prop}

An argument similar to the one we used for the Ruelle operator
$\cL_{\beta }$ in the proof of Proposition \ref{zetaRmero} now
shows

\begin{Prop}\label{zetaRmero2}
For $\beta \geqslant 0$ the Ruelle zeta function
$\zeta_{R}(z,\beta )$ for the Kac-Baker model can be written in terms of the modified
Kac-Gutzwiller operator via

$$\zeta _{R}\left( z,\beta \right)=
\prod _{\underline{\alpha }\in \left\{ 0,1\right\} ^{m}}
   \det \left( 1-z\underline{\lambda}^{\underline{\alpha}}
   \cG_{\beta }\right) ^{\left( -1\right) ^{\left|\underline{\alpha }\right| +1}}.
$$

\qed
\end{Prop}
Since the zeta function $\zeta_{R}(z,\beta)$ is meromorphic and
its divisor for fixed $\beta$ uniquely determined, it is not too
difficult to see that at least for generic values of the
parameters
 $\beta>0$ and $\underline{\lambda}\in ]0,1[^{m}$ the spectra of the two
 operators $\cL_{\beta}$ and $\cG_{\beta}$ have to be identical.
 This indeed has been shown in the case $m=1$ already by B.\,Moritz in (see \cite{Mo89}). We will
 show however (see Theorem \ref{THM}) that the spectra of the two operators
 coincide for any real $\beta$ and all parameters $\underline{\lambda}$.
 Hence the Fredholm determinants $\det(1-z\cG_{\beta})$ and
 $\det(1-z\cL_{\beta})$ of the Kac-Gutzwiller operator and the Ruelle
 operator coincide on the real axis and extend to a holomorphic function in the
 entire $\beta$ plane even if the operator $\cG_{\beta}$ contrary to the
 operator $\cL_{\beta}$ has itself no such analytic continuation to the
 entire \(\beta\)-plane.

\section{Hermite Functions and Mehler's Formula}
\label{HFMF}

Consider the operators
\[
Z_{j}=m_{x_{j}}+\frac{1}{2\pi }\frac{\partial }{\partial x_{j}}\quad,
\quad Z^{*}_{j}=m_{x_{j}}-\frac{1}{2\pi }\frac{\partial }{\partial x_{j}}\quad ,
\quad j=1,\ldots, n,\]
 where
 \( m_{x_{j}} \)
 denotes the multiplication operator
\[
\left( m_{g}f\right) \left(\underline{x}\right) =g(\underline{x})f(\underline{x})\]
 in the space
 \( L^{2}\left( \mathbb{R}^{m},d{\underline{x}}\right)  \)
for $g(\underline{x})=x_{j}$ and
 \(\underline{x}=\left( x_{1},\ldots ,x_{m}\right) \in \mathbb{R}^{m}\).

 The Hermite functions
 \( h_{\underline{\alpha }}\in L^{2}\left( \mathbb {R}^{m},d\underline{x}\right)  \)
 with
 \(\underline{\alpha }=\left( \alpha _{1},\ldots ,\alpha _{m}\right) \in \mathbb{N}^{m}_{0} \)
 are given by (see  \cite[p.51]{Fo89})
\begin{eqnarray*}
h_{\underline{0}}\left(\underline {x}\right)
&=& 2^{\frac{m}{4}}e^{-\pi \underline{x}\cdot\underline{x}}\\
h_{\underline{\alpha }}\left(\underline {x}\right)
&=& \sqrt{\frac{\pi ^{\left| \underline{\alpha }\right| }}{\underline{\alpha }!}}
    \left(\underline{Z}^{*\underline{\alpha }}h_{\underline{0}}\right) \left(\underline{x}\right)
  =  \frac{2^{\frac{m}{4}}}{\sqrt{\underline{\alpha }!}}
     \left( \frac{-1}{2\sqrt{\pi }}\right) ^{\left|\underline{\alpha }\right| }
     e^{\pi \underline{x}\cdot\underline{x}}\left( \frac{\partial }{\partial \underline{x}}\right) ^{\underline{\alpha }}
     e^{-2\pi \underline{x}\cdot \underline{x}},
\end{eqnarray*}
 where
$\underline{\alpha }=\left( \alpha _{1},\ldots ,\alpha _{m}\right)$
with  $\alpha _{i}\in \mathbb{N}_{0}$
for all  \( 1\leq i\leq m \),
$\left|\underline{\alpha }\right| =\sum ^{m}_{i=1}\alpha _{i}$ and
$\underline{\alpha }!=\prod ^{m}_{i=1}\alpha _{i}!$. Moreover,
$\underline {Z}^{*\underline{\alpha }}$ denotes the operator
$Z^{*\alpha _{1}}_{1}\cdots Z^{*\alpha _{m}}_{m}$.

The Hermite functions are known to be an orthonormal basis of the Hilbert
 space
 \( {L}^{2}\left( \mathbb{R}^{m},d\underline{x}\right)  \)
 and the following formula due to Mehler holds (see\cite{Fo89}):
\[
\sum _{\underline{ \alpha}\in \mathbb{N}_{0}^{m}} w^{\left|\underline{\alpha }\right| }
       h_{\underline{\alpha }}\left(\underline{x}\right) h_{\underline{\alpha }}
       \left(\underline{y}\right)
=\left( \frac{2}{1-w^{2}}\right) ^{\frac{m}{2}}
  \exp\left( \frac{-\pi \left( 1+w^{2}\right) \left(\underline{x}^{2}+\underline{y}^{2}\right)
            +4\pi w\,\underline{x}\cdot\underline {y}}{1-w^{2}}\right),\]
where $\left| w\right| <1$ and
$\Re\frac{2}{1-w^{2}}>0$.
>From the case
$m=1$  we then get for
$\underline{\lambda }=\left( \lambda _{1},\ldots ,\lambda _{m}\right)$
with $0<\lambda _{i}<1$
for $1\leq i\leq m$ the identity
\[
\prod ^{m}_{l=1}\sum ^{\infty }_{\alpha _{l}=0}\lambda ^{\alpha _{l}}_{l}
h_{\alpha _{l}}\left( x_{l}\right) h_{\alpha _{l}}\left( y_{l}\right)
=\prod ^{m}_{l=1}\left( \frac{2}{1-\lambda ^{2}_{l}}\right) ^{\frac{1}{2}}\cdot
\exp\left( \sum ^{m}_{l=1}\frac{-\pi \left( 1+\lambda ^{2}_{l}\right)
\left( x^{2}_{l}+y^{2}_{l}\right) +4\pi \lambda _{l}x_{l}y_{l}}{1-\lambda ^{2}_{l}}\right).\]
 A  simple calculation presented for $m=1$ in
 \cite{HiMa01}
shows that  the following proposition is true.

\begin{Prop}\label{KacKernProp}
\begin{enumerate}

\item[{\rm(i)}]
For
 $\lambda _{l}=e^{-\gamma _{l}}$ with $\Re\gamma _{l}>0$
 one has
\[
\frac{-\pi \left( 1+\lambda ^{2}_{l}\right) \left( x^{2}+y^{2}\right)
+4\pi \lambda _{l}xy}{1-\lambda ^{2}_{l}}
=\frac{2\pi }{\sinh \gamma_{l} }\left( -\frac{1}{2} \left( x^{2}+y^{2}\right)
\cosh\gamma_{l} +xy\right)
\]
\item[{\rm(ii)}]
$
\frac{1}{2\sinh \gamma_{l} }\left( -\frac{1}{2}
\left( x^{2}+y^{2}\right)\cosh\gamma_{l} +xy\right)
=-\frac{1}{4}\left(\left( x^{2}+y^{2}\right)\tanh\left(\frac{\gamma_{l}}{2}\right)
+\frac{\left( x-y\right) ^{2}}{\sinh \gamma_{l} }\right)
$
\item[{\rm(iii)}]
For
 \( \lambda _{l}=e^{-\gamma _{l}}$ and $x_{l}=\xi _{l}\frac{1}{2\sqrt{\pi }}\, ,\,
 y_{l}=\eta _{l}\frac{1}{2\sqrt{\pi }}$  one has the following version of Mehler's formula
\[
\sum _{ \underline{\alpha}\in \mathbb{N}_{0}^{m}}\underline{\lambda }^{\underline{\alpha }}h_{\underline{\alpha }}
\left( \underline{x}\right) h_{\underline{\alpha }}\left(\underline{y}\right)
=\prod ^{m}_{l=1}\left( \frac{1}{\lambda_{l}\sinh \gamma _{l}}
\right) ^{\frac{1}{2}}
\exp\left(-\frac{1}{4}\sum ^{m}_{l=1}
\left( \left( \xi ^{2}_{l}+\eta ^{2}_{l}\right)\tanh\left(\frac{\gamma_{l}}{2}\right)
       +\frac{\left( \xi _{l}-\eta _{l}\right) ^{2}}{\sinh \gamma _{l}}\right)\right)
\]
\end{enumerate}
\qed
\end{Prop}

In Proposition \ref{KacKernProp} (iii) we used only the fact that
 \[
h_{\underline{\alpha }}\left(\underline {x}\right)
=h_{\alpha _{1}}\left( x_{1}\right) \cdots h_{\alpha _{m}}\left( x_{m}\right)
\]
for
 \( \underline{\alpha }=\left( \alpha _{1},\ldots ,\alpha _{m}\right)  \),
  \( \underline{x}=\left( x_{1},\ldots ,x_{m}\right)  \)
 (see \cite[p.52]{Fo89}) and
 \(\frac{2}{1-\lambda^{2}_{l}}
 = \frac{2e^{\gamma_{l}}}{e^{\gamma_{l}}-e^{-\gamma_{l}}}
 =(\lambda_{l}\sinh{\gamma_{l}})^{-1}\).

Define next the kernel function
\begin{equation}\label{KacGutzKernmod}
\widetilde{\mathcal{K}}\left( \underline{\xi},\underline{\eta}\right)
=2\prod ^{m}_{l=1}\left( \frac{1}{4\pi\sinh \gamma _{l}}\right) ^{\frac{1}{2}}
\exp\left(-\frac{1}{4}\left( \sum ^{m}_{l=1}
      \left( \xi^{2}_{l}+\eta^{2}_{l}\right)\tanh\frac{\gamma_{l}}{2}
     +\frac{\left( \xi_{l}-\eta_{l}\right) ^{2}}{\sinh \gamma_{l}}\right) \right).
\end{equation}
Then the kernel
\(\mathcal{K}_{\beta }(\underline{\xi},\underline{\eta}) \)
of the Kac-Gutzwiller transfer operator defined in (\ref{KacGutzKern})
satisfies
\begin{equation}\label{relatekernels}
\mathcal{K}_{\beta }(\underline{\xi},\underline{\eta})=\left( \cosh \left( \sum\limits ^{m}_{l=1}
\sqrt{\beta J_{l}}\xi_{l}\right) \cosh \left( \sum\limits ^{m}_{l=1}
\sqrt{\beta J_{l}}\eta_{l}\right) \right) ^{\frac{1}{2}}\widetilde{\mathcal{K}}(\underline{\xi},\underline{\eta})
\end{equation}

\begin{Lemma}
 For an invertible real $(m\times m)$-matrix $C$
 and a smooth function
\( a:\mathbb {R}^{m}\to ]0,\infty[  \)
consider the map
\[
\left( R_{C}f\right) \left( \underline{x}\right)
=\left|\det(C)\right| ^{\frac{1}{2}}f(C\underline{x}).
\]
Then
\[
R_{C}:L^{2}\left( \mathbb {R}^{m},a(C^{-1}\underline{\xi}) d\underline{\xi }\right)
  \to L^{2}\left( \mathbb {R}^{m},a( \underline{x})
d\underline{x}\right)
\]
is an isomorphism of Hilbert spaces.
\end{Lemma}

\pf
\begin{eqnarray*}
\int _{\mathbb {R}^{m}}\left| \left( R_{C}f\right)
\left( \underline{x}\right) \right| ^{2}a\left( \underline{x}\right) d\underline{x}
&=&\int _{\mathbb {R}^{m}}\left|\det (C)\right|\, \left|
   f\left(C\underline{x}\right) \right| ^{2}\, a( \underline{x})
   d\underline{x}\\
&=&\int\limits _{\mathbb {R}^{m}}\left| f\left( \underline{\xi }\right) \right| ^{2}
   a\left(C^{-1}\underline{\xi}\right) d\underline{\xi }\\
&=&\left\Vert f\right\Vert^2 _{L^{2}\left( \mathbb {R}^{m},
   a\left(C^{-1}\underline{\xi}\right) d\underline{\xi }\right) }
\end{eqnarray*}
\qed

If
$K:L^{2}\left( \mathbb {R}^{m},d\underline{\xi }\right)
\to L^{2}\left( \mathbb {R}^{m},d\underline{\xi }\right)$
is given by an integral kernel
$K\left( \underline{\xi },\underline{\eta }\right)$
then the induced operator
$K_{C}:= R_{C}\circ K \circ R_{C}^{-1}:L^{2}\left( \mathbb {R}^{m},d\underline{x}\right) \to L^{2}
\left( \mathbb {R}^{m},d\underline{x}\right)$
has kernel
\[
K_{C}\left( \underline{x},\underline{y}\right)
=|\det(C)| K\left(C\underline{x},C\underline{y}\right)\]
as one easily verifies by a straightforward calculation using the transformation formula.

 If $\underline{c}= \left( c_1,\ldots
,c_{m}\right) \in \mathbb {R}^{m}$ with $c_{i}\neq 0$ for  $1\leq
i\leq m$ and $\boldsymbol{C}$ is the diagonal matrix with the
$c_l$ on the diagonal, we simply write $R_{\underline{c}}$ for
$R_{\boldsymbol{C}}$ and $K_{\underline{c}}$ for
$K_{\boldsymbol{C}}$. Note that
$\boldsymbol{C}\underline{x}=\left(c_{1}x_{1},\ldots,c_{m}x_{m}\right)$.

\begin{Lemma} \label{multlemma}
Let $a:\mathbb {R}^{m}\to [1,\infty[$
be a smooth function and $K$ be a bounded operator on
$L^{2}\left( \mathbb {R}^{m},d\underline{x}\right)$ given by the kernel
$K( \underline{x},\underline{y})$ as
$$\left( Kf\right) \left( \underline{x}\right)
   =\int\limits _{\mathbb {R}^{m}}K\left( \underline{x},\underline{y}\right)
    f( \underline{y}) d\underline{y}.$$
Then
\begin{enumerate}
\item[\rm{(a)}] the operator $K\circ m_{\sqrt{a}}$
                is an unbounded operator on
 $L^{2}\left( \mathbb {R}^{m},d\underline{x}\right)$
  with kernel
 $K\left( \underline{x},\underline{y}\right) \sqrt{a( \underline{y}) }$.
\item[{\rm(b)}] the operator
 $m_{\frac{1}{\sqrt{a}}}\circ K$
is a bounded operator on
$L^{2}\left( \mathbb {R}^{m},d\underline{x}\right)$
with kernel
$\frac{1}{\sqrt{a\left( \underline{x}\right) }}
K\left( \underline{x},\underline{y}\right)$.
\end{enumerate}
\end{Lemma}

\pf
\begin{enumerate}
\item[(a)]
This follows from
\[ \left( K\circ m_{\sqrt{a}}f\right) \left( \underline{x}\right)
=\int\limits _{\mathbb {R}^{m}}K\left( \underline{x},\underline{y}\right)
\sqrt{a( \underline{y}) }f( \underline{y}) d\underline{y} \]
and
\[ \int\limits _{\mathbb {R}^{m}}\left| \sqrt{a( \underline{x}) }
   f( \underline{x}) \right| ^{2}d\underline{x}
=\int\limits _{\mathbb {R}^{m}}\left| f( \underline{x}) \right| ^{2}\:
a( \underline{x}) d\underline{x}. \]

\item[(b)]
Calculate
$$
\left( m_{\frac{1}{\sqrt{a}}}\circ Kf\right) \left( x\right)
=\frac{1}{\sqrt{a( \underline{x}) }}\int\limits _{\mathbb {R}^{m}}
   K( \underline{x},\underline{y}) f( \underline{y})
   d\underline{y}
=\int\limits _{\mathbb {R}^{m}}\frac{1}{\sqrt{a(\underline{x})}}
  K( \underline{x},\underline{y}) f( \underline{y}) d\underline{y}.
$$

\end{enumerate}
\qed

Fix $a\colon \R^m\to[1,\infty[$. Then Lemma \ref{multlemma} yields
the following commutative diagram
\[
\xymatrix{
L^{2}\left( \mathbb {R}^{m},d\underline{x}\right)
&L^{2}\left( \mathbb {R}^{m},a( \underline{x}) d\underline{x}\right)
 \ar@{_(->}[l]
 \ar[r]^{K}
 \ar[d]^{m_{\sqrt{a}}}
&L^{2}\left( \mathbb {R}^{m},d\underline{x}\right)
 \ar[d]^{\id}
\\
&L^{2}\left( \mathbb {R}^{m},d\underline{x}\right)
 \ar[d]^{\id}
 \ar[r]^{K'}
& L^{2}\left( \mathbb {R}^{m},d\underline{x}\right)
 \ar[d]^{m_{\frac{1}{\sqrt{a}}}}
\\
&L^{2}\left( \mathbb {R}^{m},d\underline{x}\right)
 \ar[r]^{K''}
& L^{2}\left( \mathbb {R}^{m},d\underline{x}\right)
}
\]
with integral operators
 \( K'\)
 and
 \( K'' \)
 given by the kernels
 \[ K'\left(\underline{x},\underline{y}\right) =K\left( \underline{x},\underline{y}\right)
\frac{1}{\sqrt{a(\underline{y}) }}
\quad\mbox{ and }\quad
 K''\left(\underline{x},\underline{y}\right) =\frac{1}{\sqrt{a( \underline{x}) }}
K\left( \underline{x},\underline{y}\right) \frac{1}{\sqrt{a( \underline{y}) }}
.  \]
If the kernel $|K(x,y)|$ defines a bounded operator on $L^2(\R^m,d\underline{x})$,
then also the operator $K''$ is bounded.

Now consider the Kac-Gutzwiller operator
 \(\mathcal{K}_{\beta } \)
  with kernel
 \(\mathcal{K}_{\beta }( \underline{\xi },\underline{\eta })  \).
 For $\underline{s},\underline{x}\in \R^m$ we set
$\cosh_{\underline s}(\underline x):=\cosh(\underline{s}\cdot\underline{x})$.
With  $\underline{J}:=(J_1,\ldots,J_l)$ and
$\underline{s}=\underline{s}_0:=2\sqrt{\beta\pi} \underline{J}$
we choose the function $a:=\cosh_{\underline{s}_0}$ and for
$\underline{c}=\underline{c}_0:=2\sqrt{\pi}(1,\ldots,1)$ we obtain
\begin{eqnarray*}
\mathcal{K}_{\beta ,\underline{c}_{0}}\left( \underline{x},\underline{y}\right)
&=&\left( 2\sqrt{\pi }\right) ^{m}\mathcal{K}_{\beta }\left(2\sqrt{\pi }\underline{x},
   2\sqrt{\pi }\underline{y}\right)\\
&=&\left( 2\sqrt{\pi }\right) ^{m}\mathcal{K}_{\beta }
   \left( \underline{\xi },\underline{\eta }\right)\\
&=&2\left( \cosh_{\underline{s}_0}(\underline{x})\cosh_{\underline{s}_0}(\underline{x})
    \right)^{\frac{1}{2}}\prod ^{m}_{l=1}\left( e^{-\frac{\gamma _{l}}{2}}\right)
    \sum _{\underline{\alpha }\in\mathbb{N}_{0}^{m}}
\underline{\lambda }^{\underline{\alpha }}
  h_{\underline{\alpha }}\left( \underline{x }\right)
  h_{\underline{\alpha }}\left( \underline{y}\right)
\end{eqnarray*}

Hence for the kernel $\mathcal{K}''_{\underline{c}_{0}}(\underline{x},\underline{y})$
of the operator
$\mathcal{K}''_{\beta,\underline{c}_{0}}=:\cK''_{\underline{c}_0}$,
which does not depend on the variable $\beta$, one finds
\[\mathcal{K}''_{\underline{c}_{0}}\left(\underline{x},\underline{y}\right)
:=2\sum_{
  \underline{\alpha}\in \mathbb{N}_{0}^{m}}\prod ^{m}_{l=1}
  \left( e^{-\frac{\gamma _{l}}{2}}\right) \underline{\lambda }^{\underline{\alpha} }
   h_{\underline{\alpha }}\left( \underline{x}\right)
   h_{\underline{\alpha }}\left( \underline{y}\right).
\]
{}From the fact that the Hermite functions
 \( h_{\underline{\alpha }} \)
 determine an orthonormal basis of
 \( L^{2}\left( \mathbb {R}^{m},d\underline{x}\right)  \)
 one concludes
\begin{equation}\label{eigenfunction1}
\left(\mathcal{K}''_{\underline{c}_{0}}h_{\underline{\alpha }}\right)
\left( \underline{x}\right) =2\underline{\lambda }^{\underline{\alpha }}
\underline{\lambda }^{\underline{\frac{1}{2}}}
h_{\underline{\alpha }}( \underline{x}),
\end{equation}
i.e., the
$h_{\underline{\alpha }}$
are the complete set of eigenfunctions of the operator $\mathcal{K}''_{\underline{c}_{0}}$ with eigenvalue
\begin{equation}\label{eigenfunction2}
\rho_{\underline{\alpha}}:=
2\underline{\lambda }^{\underline{\alpha }+\underline{\frac{1}{2}}}
=2\prod ^{m}_{l=1}\lambda ^{\alpha _{l}+\frac{1}{2}}_{l}.
\end{equation}
In particular, $\mathcal{K}''_{\underline{c}_{0}}$ is bounded.
We will need this result later on.
Note that Lemma \ref{multlemma} yields the following commutative diagram for the
Kac-Gutzwiller operator $\mathcal{K}_{\beta}$:
 \[
\xymatrix{
L^{2}\left( \mathbb {R}^{m},a\left( \underline{c}_0^{-1}\circ \underline{\xi }\right) d\underline{\xi }\right)
 \ar[r]^{\quad\quad \mathcal{K}_{\beta}}
 \ar[d]_{R_{\underline{c}_0}}
& L^{2}\left( \mathbb {R}^{m}d\underline{\xi }\right)
 \ar[d]^{R_{\underline{c}_0}}
\\
L^{2}\left( \mathbb {R}^{m},a\left( \underline{x}\right) d\underline{x}\right)
 \ar[d]_{m_{\sqrt{a}}}
 \ar[r]^{\mathcal{K}_{\beta,\underline{c}_{0}}}
&L^{2}\left( \mathbb {R}^{m},d\underline{x}\right)
 \ar[d]^{\id}
\\
L^{2}\left( \mathbb {R}^{m},d\underline{x}\right)
 \ar[d]_{\id}
 \ar[r]^{\mathcal{K}'_{\beta,\underline{c}_{0}}}
&L^{2}\left( \mathbb {R}^{m},d\underline{x}\right)
 \ar[d]^{m\frac{1}{\sqrt{a}}}
\\
L^{2}\left( \mathbb {R}^{m},d\underline{x}\right)
 \ar[r]^{\mathcal{K}''_{\underline{c}_{0}}}
& L^{2}\left( \mathbb {R}^{m},d\underline{x}\right)
}\]

 The functions
 \( e^{\pi \underline{x}^{2}}h_{\underline{\alpha }}\left( \underline{x}\right)  \)
 are polynomials of degree
 \( \left| \underline{\alpha }\right|  \)
 in
 \( \underline{x} \)
 (see \cite[p.52]{Fo89}).
 In view of the estimate
\[
\cosh_{\underline{R}}(\underline{\xi})=
\sqrt{\cosh \left( \sum ^{m}_{i=1}R_{i}\xi _{i}\right) }
\leqslant e^{\frac{1}{2}\sum\limits ^{m}_{i=1}
\left| R_{i}\right| \left| \xi _{i}\right| }\]
 one concludes that the functions
 \( \underline{\xi }\mapsto
h_{\underline{\alpha }}( \underline{\xi })\sqrt{\cosh_{\underline{R}}
(\underline{\xi})}  \)
 are in
 \( L^{2}\left( \mathbb {R}^{m},d\underline{\xi }\right)  \).
 Hence all the
 \( h_{\underline{\alpha }} \)
  are contained in
 \( L^{2}\left( \mathbb {R}^{m}, \cosh_{\underline{R}}(\underline{\xi})
                d\underline{\xi }\right)  \).
This together with Mehler's formula in Proposition
\ref{KacKernProp} shows immediately that the Kac-Gutzwiller
operator $\mathcal{K}_{\beta}$ with kernel
$\mathcal{K}_{\beta}(\underline{\xi},\underline{\eta})$ in
(\ref{relatekernels}) is symmetric and positive definite for
$\beta\geqslant 0$. In fact,
\[\left(f,\mathcal{K}_{\beta}\,f\right)
=2\sum_{\underline{\alpha} \in \mathbb{N}_{0}^{m}}\prod_{l=1}^{m}
\left(\frac{\lambda_{l}}{4\pi}\right)^{\frac{1}{2}}
\underline{\lambda}^{\underline{\alpha}}\left|\int_{\mathbb{R}^{m}}
f(\underline{\xi})\sqrt{\cosh_{\underline{R}_{0}}(\underline{\xi})}
h_{\underline{\alpha}}\left(\frac{\underline{\xi}}{2\sqrt{\pi}}\right)
d\underline{\xi}\right|^{2}\geqslant0
\]
where $R_{0,i}= \sqrt{\beta J_{i}},0 \leqslant i \leqslant m$ so
that the eigenvalues of the operator $\mathcal{K}_{\beta}$ are
nonnegative for $\beta \geqslant 0$.

\section{Fock Space and Segal-Bargmann Transformation}
\label{FSSBT}

For
 \( t>0 \)
  consider the Hilbert space
 \( \mathcal{H}L^{2}\left( \mathbb {C}^{m},\mu _{t}\right)  \)
  of entire functions
 \( F:\mathbb {C}^{m}\to \mathbb {C} \)
  with
 \[
\left\Vert F\right\Vert ^{2}_{t}:=\int\limits _{\mathbb {C}^{m}}\left| F(\underline{z})\right| ^{2}
\mu_{t}(\underline{z})\, d\underline{z}<\infty, \]
 where
 \( \mu _{t} \)
  denotes the weight function
 \[
\mu _{t}\left( \underline{z}\right)
=t^{m}\exp\left( -\pi t\left| \underline{z}\right| ^{2}\right).
\]
 The Bargmann transform
 \( B_{t}:L^{2}\left( \mathbb {R}^{m},d\underline{x}\right) \to \mathcal{H}L^{2}
\left( \mathbb {C}^{m},\mu _{t}\right)  \)
defined via
\[
\left( B_{t}f\right) \left( \underline{z}\right)
=\left( \frac{2}{t}\right) ^{\frac{m}{4}}\int\limits _{\mathbb {R}^{m}}
f\left( \underline{x}\right)
\exp \left( 2\pi \underline{x}\cdot \underline{z}-\frac{\pi }{t}
\underline{x}^{2}-\frac{\pi t}{2}\underline{z}^{2}\right) d\underline{x}
\]
 determines a unitary operator in the two Hilbert-spaces (see
\cite[p.47]{Fo89}).

 In the following we are primarily interested in the case
 \( t=1 \).
 In this case we denote the space
 \( \mathcal{H}L^{2}\left( \mathbb {C}^{m},\mu _{t}\right)  \)
  simply by
 \(\mathcal{F}_{m} \)
  and call it the Fock space over
 \( \mathbb {C}^{m} \).
 The transform
 \( B_{1} \)
  is then denoted by
 \( B \)
  and hence
\begin{equation}\label{Bargmanntrafo}
\left( Bf\right) \left( \underline{z}\right)
=2^{\frac{m}{4}}\int\limits _{\mathbb {R}^{m}}f(\underline{x})
\exp \left( 2\pi \underline{x}\cdot \underline{z}-\pi \underline{x}^{2}
-\frac{\pi }{2}\underline{z}^{2}\right) d\underline{x}.
\end{equation}
Its inverse
 \( B^{-1}:\mathcal{F}_{m}\to L^{2}\left( \mathbb {R}^{m},d\underline{x}\right)  \)
  is given by
(see \cite[p.45]{Fo89})
\[
\left( B^{-1}F\right) \left( \underline{x}\right)
=2^{\frac{m}{4}}\int\limits _{\mathbb {C}^{m}}F\left( \underline{z}\right)
\exp \left( 2\pi \underline{x}\cdot \underline{z}^{*}-\pi \underline{x}^{2}-\frac{\pi }{2}
{\underline{z}^*}^{2}\right)
\exp\left( -\pi \left|\underline{z}\right| ^{2}\right)d\underline{z},
\]
where $z^*_l=\oline z_l$ simply is the complex conjugate of $z_l$.
An orthonormal basis of Fock space
 \(\mathcal{F}_{m} \)
 is given by the functions
\begin{equation}\label{ZETA}
\zeta _{\underline{\alpha }}\left( \underline{z}\right)
=\sqrt{\frac{\pi ^{\left| \underline{\alpha }\right| }}{\underline{\alpha }!}}
\underline{z}^{\underline{\alpha }},
\end{equation}
where
 \( \underline{z}^{\underline{\alpha }}=\prod\limits ^{m}_{l=1}z^{\alpha _{i}}_{i} \)
  and
\( \underline{\alpha }=\left( \alpha _{1},\ldots ,\alpha _{m}\right)
\in \mathbb {N}_{0}^{m}. \)
 Indeed one has
(see \cite[p.51]{Fo89})
\[
\zeta _{\underline{\alpha }}=Bh_{\underline{\alpha }},
\]
 where as before the
$h_{\underline{\alpha }}$
 are the Hermite functions in
 \( \mathbb {R}^{m} \).
 This we use to prove the following proposition.
\begin{Prop}\label{MlambdaProp}
 For
 \( \underline{\lambda }=\left( \lambda _{1},\ldots ,\lambda _{m}\right) \in
\left( 0,1\right) ^{m} \)
 the bounded operator
 \( M_{\underline{\lambda }}:\mathcal{F}_{m}\to\mathcal{F}_{m} \)
 defined as
 \[ M_{\underline{\lambda }}:=B\circ \mathcal{K}''_{\underline{c}_{0}}\circ B^{-1} \]
 is given by the expression
\[
\left( M_{\underline{\lambda} }F\right) \left(
\underline{z}\right) =2\sqrt{\prod ^{m}_{l=1}\lambda _{l}}F\left(
\lambda _{1}z_{1},\ldots ,\lambda _{m}z_{m} \right)
=2\underline{\lambda }^{\underline{\frac{1}{2}}} F\left(
\boldsymbol{\Lambda }\underline{z}\right),
 \]
 where
 \( \underline{\frac{1}{2}}=\left( \frac{1}{2},\frac{1}{2},\ldots, \frac{1}{2}\right)  \)
 and
 \( \boldsymbol{\Lambda } \underline{z}
:=\left( \lambda _{1}z_{1},\ldots, \lambda _{m}z_{m}\right)  \).
\end{Prop}

\pf
Consider first
$F( \underline{z}) =\zeta _{\underline{\alpha }}( \underline{z})$.
In view of (\ref{eigenfunction1}) and (\ref{eigenfunction2}) we have
 \[
M_{\underline{\lambda }}\zeta _{\underline{\alpha }}
=B\circ\mathcal{K}''_{\underline{c}_{0}}\circ B^{-1}\zeta _{\underline{\alpha }}
=B\circ\mathcal{K}''_{\underline{c}_{0}}h_{\underline{\alpha} }=B\left( 2\underline{\lambda }^{\underline{\alpha }
+\underline{\frac{1}{2}}}h_{\underline{\alpha }}\right)
=2\underline{\lambda }^{\underline{\alpha }+\underline{\frac{1}{2}}}
\zeta _{\underline{\alpha }}.\]
 But
 \( \zeta _{\underline{\alpha} }\left( \underline{z}\right)  \)
 is homogeneous of degree
 \( \alpha _{i} \)
 in
 \( z_{i} \)
 and hence
\[
\underline{\lambda}^{\underline{\alpha }} \zeta
_{\underline{\alpha }}\left( \underline{z}\right) = \zeta
_{\underline{\alpha }}\left(
\boldsymbol{\Lambda}\underline{z}\right).\]
 Therefore the claim is true for the basis elements
 \( \zeta _{\underline{\alpha }} \)
 of
 \(\mathcal{F}_{m} \)
 and hence also for any
 \( F\in \mathcal{F}_{m} \).
\qed
 For
 \( \underline{r}\in \mathbb {R}^{m} \)
 define the translation operator
$\tau _{\underline{r}}:L^{2}\left( \mathbb {R}^{m},d\underline{x}\right)
\to L^{2}\left( \mathbb {R}^{m},dx\right)$ by
$$
\left( \tau _{\underline{r}}f\right) \left( \underline{x}\right)
:=f\left( \underline{x}-\underline{r}\right)
$$
and for
$s\in\R\setminus\{0\}$
define the multiplication operator
$\mu _{s}:L^{2}\left( \mathbb {R}^{m},d\underline{x}\right)
\to L^{2}\left( \mathbb {R}^{m},d\underline{x}\right)
$
by
$$\left( \mu _{s}f\right) \left( \underline{x}\right)
:=sf\left( \underline{x}\right)
$$
For $\underline{s}\in(\R\setminus\{0\})^m$ we define
$\mu_{\underline{s}}^{\underline{\alpha}}:=\prod_{l=1}^m \mu_{s_l}^{\alpha_l}$.
Since $Z^*_l$ and $\mu_s$ commute, it makes sense to write
$\left( \underline{Z}^{*}+\mu _{\underline{s}}\right) ^{\underline{\alpha }}
:=\prod\limits ^{m}_{j=1}\left( Z^{*}_{j}+\mu_{ s_{j}}\right) ^{\alpha _{j}}$.

\begin{Prop}\label{transprop} For any \(\underline{\alpha}\in\mathbb{N}_{0}^{m}\) we have
$$
 {\underline{Z}^*}^{\underline{\alpha }}\circ \tau _{\underline{r}}
=\tau _{\underline{r}}\circ \left( \underline{Z}^{*}
+\mu _{\underline{r}}\right) ^{\underline{\alpha }},
$$

\end{Prop}

\pf
For
 \( f\in \mathcal{C}^{1}\left( \mathbb {R}^{m}\right)  \)
  and
 \( 1\leq j\leq m \)
we calculate
\begin{eqnarray*}
\left( Z^{*}_{j}\circ \tau _{\underline{r}}f\right) \left( \underline{x}\right)
&=&\left( x_{j}-\frac{1}{2\pi }\frac{\partial }{\partial x_{j}}\right)
    f( \underline{x}-\underline{r})\\
&=&\left( x_{j}-r_{j}\right) f( \underline{x}-\underline{r}) -\frac{1}{2\pi }
   \frac{\partial }{\partial x_{j}}f( \underline{x}-\underline{r})
    +r_{j}f( \underline{x}-\underline{r})\\
&=&\left( \tau _{\underline{r}}\circ Z^{*}_{j}f\right) \left( \underline{x}\right)
   +\left( \tau _{\underline{r}}\circ \mu _{r_{j}}f\right)
   \left( \underline{x}\right)\\
&=&\left( \tau _{\underline{r}}\circ \left( Z^{*}_{j}+\mu _{r_{j}}\right) f\right)
   \left( \underline{x}\right)
\end{eqnarray*}
{}From this it follows immediately that
 \( {\underline{Z}^*}^{\underline{\alpha }}\circ \tau _{\underline{r}}
=\tau _{\underline{r}}\circ \left( \underline{Z}^{*}
+\mu _{\underline{r}}\right) ^{\underline{\alpha }}\).
\qed

For the Hermite function
 \( h_{\underline{0}}\left( \underline{x}\right)  \)
one now  gets
\[ \left( {\underline{Z}^*}^{\underline{\alpha }}\circ
\tau _{\underline{r}}h_{\underline{0}}\right) \left( \underline{x}\right)
=\tau _{\underline{r}}\circ \left( \underline{Z}^{*}
+\mu _{\underline{r}}\right) ^{\underline{\alpha }}
h_{\underline{0}}\left( \underline{x}\right)
=\tau _{\underline{r}}\sum\limits ^{\underline{\alpha }}_{\underline{l}=\underline{0}}
\binom{\underline{\alpha }}{\underline{l}}\mu ^{\underline{l}}_{\underline{r}}
{\underline{Z}^*}^{(\underline{\alpha }-\underline{l})}\underline{h}_{0}
\left( \underline{x}\right),  \]
where we used the notation
\[ \sum\limits ^{\underline{\alpha }}_{\underline{l}=\underline{0}}
\binom{\underline{\alpha }}{\underline{l}}\mu ^{\underline{l}}_{\underline{r}}
{\underline{Z}^*}^{(\underline{\alpha }-\underline{l})}=\sum\limits ^{\alpha _{1}}_{l_{1}=0}
\cdots \sum\limits ^{\alpha _{m}}_{l_{m}=0}\binom{\alpha _{1}}{l_{1}}\cdots
\binom{\alpha _{m}}{l_{m}}\mu ^{l_{1}}_{r_{1}}\cdots \mu ^{l_{m}}_{r_{m}}
{Z^*_{1}}^{(\alpha _{1}-l_{1})}\cdots {Z^*_{m}}^{(\alpha _{m}-l_{m})}.
\]
For
 \( {\underline{Z}^*}^{(\underline{\alpha }-\underline{l})}h_{\underline{0}}
=:\widetilde{h}_{\underline{\alpha }-\underline{l}} \)
we find
\[ {\underline{Z}^*}^{\underline{\alpha }}\circ \tau _{\underline{r}}h_{\underline{0}}
=\sum\limits ^{\underline{\alpha }}_{\underline{l}=\underline{0}}
\binom{\underline{\alpha }}{\underline{l}}\mu ^{\underline{l}}_{\underline{r}}
\tau _{\underline{r}}\widetilde{h}_{\underline{\alpha }-\underline{l}}.
\]
For \(\underline{s}\in \mathbb{R}^{m}\) denote by
$\exp_{\underline{s}}\colon \R^m\to\R$ the function defined by
\(\exp_{\underline{s}}(\underline{x}):= e^{\underline{s}\cdot\underline{x}}\). Then one has
\begin{Prop}\label{emultprop}
 For
 \( \underline{\alpha }\in \mathbb {N}_{0}^{m} \)
  and
 \( \underline{s}\in \mathbb {R}^{m} \)
  the following identities hold
\begin{eqnarray*}
\underline{Z}^{*\underline{\alpha }}\circ m_{\exp_{\underline{s}}}
&=&m_{\exp_{\underline{s}}}
\circ  \left( \underline{Z}^{*}-\frac{\underline{s}}{2\pi }\right)^{\underline{\alpha }}
\\
 m_{\exp_{\underline{s}}}\circ \underline{Z}^{*\underline{\alpha }}
&=&\left( \underline{Z}^{*}+\frac{\underline{s}}{2\pi }\right) ^{\underline{\alpha }}\circ
m_{\exp_{\underline{s}}}
\end{eqnarray*}
\end{Prop}

 \pf
By definition of
 \( Z^{*}_{j} \)
  we obtain for smooth
 \( f \)
\begin{eqnarray*}
\left( Z^{*}_{j}\circ m_{\exp_{\underline{s}}} f\right)
\left( \underline{x}\right)
&=&x_{j}e^{\underline{s}\cdot \underline{x}}
   f\left( \underline{x}\right) -\frac{1}{2\pi }\frac{\partial }{\partial x_{j}}
   \left( e^{\underline{s}\cdot \underline{x}}f\left( \underline{x}\right) \right)\\
&=&x_{j}e^{\underline{s}\cdot \underline{x}}f\left( \underline{x}\right)
   -\frac{1}{2\pi }
   \left( s_{j}e^{\underline{s}\cdot \underline{x}}f\left( \underline{x}\right)
    +e^{\underline{s}\cdot \underline{x}}\frac{\partial }{\partial x_{j}}
    f\left( \underline{x}\right) \right)\\
&=&e^{\underline{s}\cdot \underline{x}}\left( Z^{*}_{j}f\right)
    \left( \underline{x}\right) -\frac{s_{j}}{2\pi }
   e^{\underline{s}\cdot \underline{x}}f\left( \underline{x}\right)\\
&=&\left(m_{\exp_{\underline{s}}} \circ\left( Z^{*}_{j}-\frac{s_{j}}{2\pi }\right)
   f\right)\left( \underline{x}\right).
\end{eqnarray*}
 Iterating this calculation proves the first identity of the proposition.
 The second identity is proved in the same way.
\qed

{}From this one derives

\begin{Prop}
For
 \( \underline{s}\in \mathbb {R}^{m} \)
 and
 \( \widetilde{h}_{\underline{\alpha} }={\underline{Z}^*}^{\underline{\alpha} }h_{\underline{0}} \)
 one has
\[ m_{\exp_{\underline{s}}} \widetilde{h}_{\underline{\alpha }}
  =e^{\frac{\underline{s}^{2}}{4\pi }}\sum\limits ^{\underline{\alpha }}_{\underline{k}
  =\underline{0}}\binom{\underline{\alpha }}{\underline{k}}\left( \frac{\underline{s}}{\pi }
  \right) ^{\underline{\alpha }-\underline{k}}\tau _{\frac{\underline{s}}{2\pi }}
   \widetilde{h}_{\underline{k}}.
\]
\end{Prop}

\pf
For
 \( \underline{\alpha }=\underline{0} \)
  we have
 \( \widetilde{h}_{\underline{0}}=h_{\underline{0}} \)
and hence we get
\begin{eqnarray*}
\left(m_{\exp_{\underline{s}}} h_{\underline{0}}\right)\left(\underline{x}\right)
&=&e^{\underline{s}\cdot \underline{x}} 2^{\frac{m}{4}}e^{-\pi \underline{x}^{2}}\\
&=&2^{\frac{m}{4}}e^{\underline{s}\cdot \underline{x}-\pi \underline{x}^{2}}\\
&=&2^{\frac{m}{4}}e^{-\pi \left( \underline{x}-\frac{\underline{s}}{2\pi }\right) ^{2}
+ \frac{\underline{s}^{2}}{4\pi }}\\
&=&e^{\frac{\underline{s}^{2}}{4\pi } }\tau _{\frac{\underline{s}}{2\pi }}
    h_{\underline{0}}\left( \underline{x}\right).
\end{eqnarray*} For general
 \( \underline{\alpha }\in \mathbb {N}^{m}_{0} \)
  we then get using Proposition \ref{emultprop} and Proposition \ref{transprop}
\begin{eqnarray*}
 m_{\exp_{\underline{s}}}\widetilde{h}_{\underline{\alpha }}
&=&m_{\exp_{\underline{s}}}\circ {\underline{Z}^*}^{\underline{\alpha} }
   h_{\underline{0}}\\
&=&\left( \underline{Z}^{*}+\frac{\underline{s}}{2\pi }\right) ^{\underline{\alpha }}
   m_{\exp_{\underline{s}}}h_{\underline{0}}\\
&=&e^{\frac{\underline{s}^{2}}{4\pi } }
  \left( \underline{Z}^{*}
  +\frac{\underline{s}}{2\pi }\right) ^{\underline{\alpha }}
  \circ \tau _{\frac{\underline{s}}{2\pi }}h_{\underline{0}}\\
&=&e^{\frac{\underline{s}^{2}}{4\pi }}
   \sum\limits ^{\underline{\alpha }}_{\underline{l}=\underline{0}}
   \binom{\underline{\alpha }}{\underline{l}}\left( \frac{\underline{s}}{2\pi }
   \right) ^{\underline{\alpha }-\underline{l}}{\underline{Z}^*}^{\underline{l}}
   \circ \tau _{\frac{\underline{s}}{2\pi }}h_{\underline{0}}\\
&=&e^{\frac{\underline{s}^{2}}{4\pi } }
   \sum\limits ^{\underline{\alpha }}_{\underline{l}=\underline{0}}
   \binom{\underline{\alpha }}{\underline{l}}\left( \frac{\underline{s}}{2\pi }
   \right) ^{\underline{\alpha }-\underline{l}}\tau _{\frac{\underline{s}}{2\pi }}
   \circ\left( \underline{Z}^{*}
   +\mu _{\frac{\underline{s}}{2\pi }}\right) ^{\underline{l}}h_{\underline{0}}\\
&=&e^{\frac{\underline{s}^{2}}{4\pi } }
   \sum\limits ^{\underline{\alpha }}_{\underline{l}=\underline{0}}
   \binom{\underline{\alpha }}{\underline{l}}
   \left( \frac{\underline{s}}{2\pi }\right) ^{\underline{\alpha }-\underline{l}}
   \sum\limits ^{\underline{l}}_{\underline{k}=\underline{0}}
   \binom{\underline{l}}{\underline{k}}\left( \frac{\underline{s}}{2\pi }
   \right) ^{\underline{l}-\underline{k}}\tau _{\frac{\underline{s}}{2\pi }}\circ
    {\underline{Z}^*}^{\underline{k}}h_{\underline{0}}\\
&=&e^{\frac{\underline{s}^{2}}{4\pi } }
   \sum\limits ^{\underline{\alpha }}_{\underline{l}=\underline{0}}
   \sum\limits ^{\underline{l}}_{\underline{k}=\underline{0}}
   \binom{\underline{\alpha }}{\underline{l}}\binom{\underline{l}}{\underline{k}}
   \left( \frac{\underline{s}}{2\pi }\right) ^{\underline{\alpha }-\underline{k}}
   \tau _{\frac{\underline{s}}{2\pi }}\widetilde{h}_{\underline{k}}\\
&=&e^{\frac{\underline{s}^{2}}{4\pi } }
   \sum\limits ^{\underline{\alpha }}_{\underline{k}=\underline{0}}
   \sum\limits ^{\underline{\alpha }}_{\underline{l}=\underline{k}}
   \binom{\underline{\alpha }}{\underline{l}}\binom{\underline{l}}{\underline{k}}
   \left( \frac{\underline{s}}{2\pi }\right) ^{\underline{\alpha }-\underline{k}}
   \tau _{\frac{\underline{s}}{2\pi }}\widetilde{h}_{\underline{k}}.
\end{eqnarray*}
But
 \( \sum\limits ^{\underline{\alpha }}_{\underline{l}=\underline{k}}
\binom{\underline{\alpha }}{\underline{l}}
\binom{\underline{l}}{\underline{k}}=\underline{2}^{\underline{\alpha }-\underline{k}}
\binom{\underline{\alpha }}{\underline{k}} \)
since
 \( \sum\limits ^{\underline{\alpha }_{i}}_{\underline{l}_{i}=\underline{k}_{i}}
\binom{\underline{\alpha }_{i}}{\underline{l}_{i}}
\binom{\underline{l}_{i}}{\underline{k}_{i}}=\underline{2}^{\underline{\alpha }_{i}
-\underline{k}_{i}}\binom{\underline{\alpha }_{i}}{\underline{k}_{i}} \),
and therefore
 \[ m_{\exp_{\underline{s}}}\widetilde{h}_{\alpha }
=e^{\frac{\underline{s}^{2}}{4\pi } }
\sum\limits ^{\underline{\alpha }}_{\underline{k}=\underline{0}}
\underline{2}^{\underline{\alpha }-\underline{k}}\binom{\underline{\alpha }}{\underline{k}}
\left( \frac{\underline{s}}{2\pi }\right) ^{\underline{\alpha }-\underline{k}}
\tau _{\frac{\underline{s}}{2\pi }}\widetilde{h}_{\underline{k}}
=e^{\frac{\underline{s}^{2}}{4\pi } }
\sum ^{\underline{\alpha }}_{\underline{k}=\underline{0}}
\binom{\underline{\alpha }}{\underline{k}}
\left( \frac{\underline{s}}{\pi }\right) ^{\underline{\alpha }-\underline{k}}
\tau _{\frac{\underline{s}}{2\pi }}\widetilde{h}_{\underline{k}}.
\]

\qed

\begin{Prop}
 For
 \( \underline{r}\in \mathbb {R}^{m} \)
 and
 \( B:L^{2}(\mathbb{R}^{m},d\underline{x}) \to \mathcal{F}_{m} \)
  the Bargmann transform one has
 \[
B\circ \tau _{\underline{r}}
=e^{-\frac{\pi }{2}\underline{r}^{2}} m_{\exp_{\pi\underline{r}}}
\circ \tau _{\underline{r}}\circ B,\]
where we have denoted the function \(e^{\underline{r}\cdot\underline{z}}\) also by \(\exp_{\underline{r}}.\) \end{Prop}

\pf
Using (\ref{Bargmanntrafo}) we calculate
\begin{eqnarray*}
\left( B\circ \tau _{\underline{r}}f\right) \left( \underline{z}\right)
&=&2^{\frac{m}{4}}\int\limits _{\mathbb {R}^{m}}
   f\left( \underline{x}-\underline{r}\right)
   e^{2\pi \underline{x}\cdot \underline{z}-\pi \underline{x}^{2}-
   \frac{\pi }{2}\underline{z}^{2}}d\underline{x}  \\
&=&2^{\frac{m}{4}}\int\limits _{\mathbb {R}^{m}}f(\underline{y})
   e^{2\pi \left( \underline{y}+\underline{r}\right)
   \cdot \underline{z}-\pi \left( \underline{y}+\underline{r}\right) ^{2}
   -\frac{\pi }{2}\underline{z}^{2}}d\underline{y}  \\
&=& 2^{\frac{m}{4}}\int\limits _{\mathbb {R}^{m}}f\left( \underline{y}\right)
    e^{2\pi \underline{y}\cdot( \underline{z}-\underline{r})
       -\pi \underline{y}^{2}-\frac{\pi }{2}( \underline{z}-\underline{r})^{2}}
    e^{\pi \underline{r}\cdot \underline{z}-\frac{\pi }{2}\underline{r}^{2}} d\underline{y}\\
&=&e^{\pi \underline{r}\cdot \underline{z}-\frac{\pi }{2}\underline{r}^{2}}
   \left( \tau _{\underline{r}}\circ Bf\right) \left(\underline{z}\right)\\
&=&e^{-\frac{\pi\underline{r}^{2}}{2}}\left(m_{\exp_{\pi\underline{r}}}\circ \tau_{\underline{r}}\circ B\,f\right)(\underline{z}).
\end{eqnarray*}
\qed

 Consider next the operator in
 \(\mathcal{F}_{m} \) induced from
 the multiplication operator
 \( m_{\exp_{\underline{s}}} \)
 in
 \( L^{2}\left( \mathbb {R}^{m},d\underline{x}\right)  \).
 One finds

\begin{Prop}\label{BSintertwine}
 For
 \( \underline{s}\in \mathbb {R}^{m} \)
 and $F\colon \C^m\to\C$ polynomial we have
 \[\left(
B\circ m_{\exp_{ \underline{s}}}\circ B^{-1}F\right)(\underline{z})
=e^{\frac{\underline{s}^{2}}{8\pi }}
 e^{\frac{\underline{s}\cdot \underline{z}}{2}}
 F\left( \underline{z}+\frac{\underline{s}}{2\pi }\right) \]

\end{Prop}

\pf
For the functions
\(\widetilde{\zeta}_{\underline{\alpha}}\left(\underline{z}\right)
:= \underline{z}^{\underline{\alpha}}
=B\tilde h_{\underline{\alpha}}(\underline{z})\) one finds
\begin{eqnarray*}
B\circ m_{\exp_{\underline{s}}} \circ B^{-1} \widetilde{\zeta}_{\underline{\alpha}}
&=&B\circ m_{\exp_{\underline{s}}} \widetilde{h}_{\underline{\alpha }}\\
&=&B e^{\frac{\underline{s}^{2}}{4\pi }}
   \sum\limits ^{\underline{\alpha }}_{\underline{k}=\underline{0}}
   \binom{\underline{\alpha }}{\underline{k}}
   \left( \frac{\underline{s}}{\pi }\right) ^{\underline{\alpha }-\underline{k}}
   \tau _{\frac{\underline{s}}{2\pi }}\widetilde{h}_{\underline{k}}\\
&=&e^{\frac{\underline{s}^{2}}{4\pi } }
   \sum\limits ^{\underline{\alpha }}_{\underline{k}=\underline{0}}
   \binom{\underline{\alpha }}{\underline{k}}
   \left( \frac{\underline{s}}{\pi }\right) ^{\underline{\alpha }-\underline{k}}
   e^{-\frac{\underline{s}^{2}}{8\pi}} m_{\exp_{\frac{\underline{s}}{2}}}
   \circ \tau _{\frac{\underline{s}}{2\pi }}\circ B\widetilde{h}_{\underline{k}} \\
&=&e^{\frac{\underline{s}^{2}}{8\pi }}e^{\frac{\underline{s}\cdot \underline{z}}
   {2}}\sum\limits ^{\underline{\alpha }}_{\underline{k}=\underline{0}}
   \binom{\underline{\alpha }}{\underline{k}}
   \left( \frac{\underline{s}}{\pi }\right) ^{\underline{\alpha }-\underline{k}}
   \left( \underline{z}-\frac{\underline{s}}{2\pi }\right) ^{\underline{k }} \\
&=&e^{\frac{\underline{s}^{2}}{8\pi }+\frac{\underline{s}\cdot \underline{z}}{2}}
   \left( \underline{z} +\frac{\underline{s}}{2\pi }\right)^{\underline{\alpha }}.
\end{eqnarray*}
Since the
$\left\{\widetilde{\zeta}_{\underline{\alpha }}\right\}$
form a basis in the space
$\cF_{m}$ the claim of the proposition is true.
\qed

\begin{Rem} According to Proposition \ref{BSintertwine}
we can view
$B\circ m_{\exp_{\underline{s}}}\circ B^{-1}$
as an unbounded operator on $\cF_m$ which is defined on a dense linear subspace, namely the
space of polynomial functions.
\qed
\end{Rem}

For the following we need the densely defined unbounded operator
$C_{\underline{s}}:\cF_{m}\to \cF_{m}$
defined as
\begin{equation}\label{Cdef}
C_{\underline{s}}:=B\circ m_{\cosh_{ \underline{s}}}\circ B^{-1}
=\frac{1}{2}\left( B\circ m_{\exp_{\underline{s}}}\circ B^{-1}
 +B\circ m_{\exp_{-\underline{s}}}\circ B^{-1}\right).
\end{equation}

>From Proposition \ref{BSintertwine}
we deduce
\[
\left( C_{\underline{s}}F\right) \left( \underline{z}\right)
=\frac{1}{2}e^{ \frac{\underline{s^{2}}}{8\pi }}
\left( e^{\frac{\underline{s}\cdot \underline{z}}{2}}
F\left( \underline{z}+\frac{\underline{s}}{2\pi }\right)+e^{-\frac{\underline{s}\cdot \underline{z}}{2}}
F\left( \underline{z}-\frac{\underline{s}}{2\pi }\right) \right)
\]
for polynomial $F$, so indeed $C_{\underline{s}}$ is densely defined.
Composing this unbounded operator with the bounded operator
$M_{\underline{\lambda }}$
of Proposition \ref{MlambdaProp} we actually get a bounded operator on
$\cF_{m}$ as  the following proposition shows.

\begin{Prop}
For
$\underline{s}\in \mathbb {R}^{m}$,
$\underline{\lambda }
=\left( \lambda _{1},\ldots ,\lambda _{m}\right) \in \left( 0,1\right) ^{m}$,
and polynomial $F\in\cF_m$ one has
\[
\left( C_{\underline{s}}\circ M_{\underline{\lambda }}F\right)
\left( \underline{z}\right)
=\underline{\lambda}^{\underline{\frac{1}{2}}}e^{\frac{\underline{s}^{2}}{8\pi }}
  \left( e^{\frac{\underline{s}\cdot \underline{z}}{2}}
  F\left( \boldsymbol{\Lambda }\underline{z}+\frac{\boldsymbol{\Lambda }
  \underline{s}}{2\pi }\right) + e^{ -\frac{\underline{s}\cdot \underline{z}}{2}}
  F\left( \boldsymbol{\Lambda }\underline{z}-\frac{\boldsymbol{\Lambda }
  \underline{s}}{2\pi }\right) \right)
\]
and hence
$C_{\underline{s}}\circ M_{\underline{\lambda }}$ extends
to a bounded operator $\cF_{m}\to \cF_{m}$.
\end{Prop}

\pf
\[
\left( \underline{C}_{\underline{s}}\circ M_{\underline{\lambda }}F\right)
  \left( \underline{z}\right)
=\frac{1}{2}e^{ \frac{\underline{s}^{2}}{8\pi }}\left( e^{ \frac{\underline{s}\cdot \underline{z}}{2}}
 \left( M_{\underline{\lambda }}F\right)
 \left( \underline{z}+\frac{\underline{s}}{2\pi }\right)
 +e^{ -\frac{\underline{s}\cdot \underline{z}}{2}}\left( M_{\lambda }F\right)
 \left( \underline{z}-\frac{\underline{s}}{2\pi }\right) \right).
\]
But we have $\left( M_{\underline{\lambda }}F\right) \left(
\underline{z}\right)
=2\underline{\lambda}^{\underline{\frac{1}{2}}}F\left(
\boldsymbol{\Lambda } \underline{z}\right)$ and hence the claim
follows. \qed
 For $\underline{\alpha }\in \mathbb {R}_{*}^{m}$ with
$\mathbb {R}_{*}=\left\{r\in \mathbb {R}: r \neq 0\right\}$ define
the map $\nu_{\underline{\alpha
}}:\cF_{m}^{\underline{\alpha}^{-1}}\to \cF_{m}$ by
\begin{equation}\label{nudef}
\left( \nu_{\underline{\alpha }}F\right) \left(
\underline{z}\right) :=F\left( \boldsymbol{A}\underline{z}\right)
=F\left( \alpha _{1}z_{1},\ldots ,\alpha _{m}z_{m}\right),
\end{equation}
where
\[
\cF^{\underline{\alpha }}_{m}
=\left\{ F:\mathbb {C}^{m}\to \mathbb {C}\quad \textrm{entire },
 \quad \left\Vert F\right\Vert ^{2}_{\underline{\alpha }}:=
 \int _{\mathbb {C}^{m}}\left| F(z)\right| ^{2}
  e^{-(\pi \sum ^{m}_{i=1}\left| \alpha _{i}z_{i}\right| ^{2})}dz<\infty \right\}
\]
and $\boldsymbol{A}$ is the diagonal matrix with diagonal entries
$\alpha_1,\ldots,\alpha_m$.

 Then consider the induced operator
\[
\nu_{\underline{\alpha }^{-1}}\circ C_{\underline{s}}\circ M_{\underline{\lambda }}
\circ \nu_{\underline{\alpha }}:
\cF^{\underline{\alpha }^{-1}}_{m}\to \cF^{\underline{\alpha }^{-1}}_{m}\quad .
\]
Inserting the expressions for  $\nu_{\underline{\alpha }}$
and $C_{\underline{s}}\circ M_{\underline{\lambda }}$
one gets

\begin{eqnarray*}
\left( \nu_{\underline{\alpha }^{-1}}\circ C_{\underline{s}}\circ
   M_{\underline{\lambda }}\circ \nu_{\underline{\alpha }}F\right)
   \left( \underline{z}\right)
&=&\underline{\lambda }^{\underline{\frac{1}{2}}} e^{
\frac{\underline{s}^{2}}{8\pi}}
  \left( e^{ \frac{\underline{s}\cdot
  \left(\boldsymbol{A}^{-1} \underline{z}\right) }{2}}
  F\left( \boldsymbol{\Lambda }\underline{z}
  +\frac{\boldsymbol{\Lambda }\boldsymbol{A}\underline{s}}{2\pi }\right)
  +e^{ -\frac{\underline{s}\cdot
  \left( \boldsymbol{A}^{-1}\underline{z}\right) }{2}}
  F\left( \boldsymbol{\Lambda }\underline{z}-\frac{\boldsymbol{\Lambda }
  \boldsymbol{A}\underline{s}}{2\pi }\right) \right).
\end{eqnarray*}
Choose next the parameters $\underline{s}=\underline{s}_{0}$ and
$\underline{\alpha }=\underline{\alpha}_{0}$ with
\[
s_{0,i}=2\sqrt{\pi \beta J_{i}}\quad \textrm{and }\quad \alpha_{0,i}
=\sqrt{\frac{\pi }{\beta J_{i}}},\quad 1\leq i\leq m.
\]
We then get
\begin{eqnarray*}
\left( \nu_{\underline{\alpha }_{0}^{-1}}\circ
C_{\underline{s}_{0}}\circ
 M_{\underline{\lambda }}\circ \nu_{\underline{\alpha _{0}}}F\right)
 \left( \underline{z}\right) =
\prod ^{m}_{i=1}\left( \lambda _{i}\exp \beta J_{i}\right) ^{\frac{1}{2}}
  \left( e^{ \beta \sum ^{m}_{i=1} J_{i}z_{i}}F
  \left( \boldsymbol{\Lambda }\underline{z}+\underline{\lambda }\right)
  +e^{-\beta \sum ^{m}_{i=1}J_{i}z_{i}}F
  \left( \boldsymbol{\Lambda }\underline{z}-\underline{\lambda }\right)\right) .
\end{eqnarray*}
But this operator has up to the multiplicative factor
$\prod\limits ^{m}_{i=1}\left( \lambda _{i}\exp \beta J_{i}\right)
^{\frac{1}{2}}$ exactly the form of the Ruelle transfer operator
of the Kac model for the parameters \(\underline{\lambda}\).

The operator
$\nu_{\underline{\alpha }_{0}^{-1}}\circ C_{\underline{s}_{0}}\circ M_{\underline{\lambda }}
\circ \nu_{\underline{\alpha} _{0}}$
is defined in the Hilbert space \(\cF^{\underline{\alpha }_{0}^{-1}}_{m}\)
with \(\alpha_{0,i}^{-1}=\sqrt{\frac{\beta J_{i}}{\pi}}\).
All eigenfunctions of the operator
$\cL_{\beta }:\mathcal{B}(D)\to \mathcal{B}(D)$
besides the ones belonging to the eigenvalue zero belong to this space.

This can be seen as follows. From the functional equation

 $$\rho f(\underline{z}) =e^{ \beta
\underline{J}\cdot\underline{z}}
f(\boldsymbol{\Lambda}\underline{z}+\underline{\lambda}) + e^{
-\beta \underline{J}\cdot\underline{z}}
f(\boldsymbol{\Lambda}\underline{z}-\underline{\lambda}) $$

 one
concludes that for $\rho \neq 0 $ any eigenfunction of
\(\cL_{\beta}\) is an entire function in $\underline{z}$ which can
grow for $\left| \underline{z}\right| \to \infty$ at most like
$e^{C \sum\limits ^{m}_{l=1}\left| z_{i}\right|}$ for some
positive constant $C$. Such functions, however, belong to any of
the Fock spaces $\cF^{\underline{\alpha }}_{m}$. Hence the
operators $\cL_{\beta }$ and $\frac{1}{ \prod\limits
^{m}_{i=1}\left( \lambda _{i}\exp \beta J_{i}\right)
^{\frac{1}{2}}} \nu_{\underline{\alpha }_{0}^{-1}}\circ
C_{\underline{s}_{0}}\circ M_{\underline{\lambda }}\circ
\nu_{\underline{\alpha }_{0}}$ have the same spectra in this
space. The eigenfunctions with eigenvalue zero of the operator
$\cL_{\beta}:\mathcal{B}(D)\to \mathcal{B}(D)$ can be determined
explicitly. They are given by the functions
\[f_{\underline{n},\underline{\alpha}}(\underline{z})= e^{-\left(\beta\sum^{m}_{l=1}
\frac{J_{l}}{2 \lambda_{l}^{2}} z_{l}^{2}\right)} \prod^{m}_{l=1}
\left(\exp\frac{(2n_{l}+1)\pi i z_{l}}{2\lambda_{l}}\right)^{\alpha_{l}}
\]
 with $\underline{n}\in \mathbb{N}_{0}^{m}$
and $\underline{\alpha} \in \mathbb{Z}^{m}$ such that $\left|\underline{\alpha}\right|=1\mod 2 $
and hence do not belong to the space $\mathcal{F}_{m}^{\underline{\alpha}_{0}^{-1}}.$
Summarizing we have shown

\begin{Prop} \label{OPCM}

The operators $\prod_{l=1}^{m}\left(\lambda_{l} \exp{\beta J_{l}}\right)^{\frac{1}{2}} 
\cL_{\beta}:\mathcal{F}_{m}^{\underline{\alpha}^{-1}_{0}} \to
\mathcal{F}_{m}^{\underline{\alpha}^{-1}_{0}} $ and $C_{\underline{s}_{0}}\circ M_{\underline{\lambda}}:
\mathcal{F}_{m}\to\mathcal{F}_{m} $ are conjugate.
\qed
\end{Prop}

This leads to

\begin{Prop}\label{Lconjugate}
The operators
$\cL_{\beta }:\cF_{m}^{\underline{\alpha}_{0}^{-1}} \to
 \cF^{\underline{\alpha}_{0}^{-1}}_{m}$
and
\[\frac{1}{\prod^{m}_{l=1}\left(\lambda_{l}\exp \beta J_{l}\right)^{\frac{1}{2}}}
m_{\cosh_{\underline{s}_{0}}}
   \circ \mathcal{K}''_{\underline{c}_{0}}:
   L^{2}\left( \mathbb {R}^{m},d\underline{x}\right)
   \to L^{2}\left( \mathbb {R}^{m},d\underline{x}\right)
\]
are conjugate.
\end{Prop}

\pf
Inserting the definitions of the operator
$C_{\underline{s}}$
and $M_{\underline{\lambda }}$
we find
\[
C_{\underline{s}}\circ M_{\underline{\lambda }}
=B\circ m_{\cosh_{ \underline{s}}}\circ B^{-1}\circ B
 \circ \mathcal{K}''_{\underline{c}_{0}}\circ B^{-1}
=B\circ m_{\cosh_{ \underline{s}}}\circ \mathcal{K}''_{\underline{c}_{0}}
  \circ B^{-1}
\]
and hence
$m_{\cosh_{\underline{s}}}\circ\mathcal{K}''_{\underline{c}_{0}}$
is conjugate to
$C_{\underline{s}}\circ M_{\underline{\lambda }}$. Therefore
$f\in L^{2}\left( \mathbb {R}^{m},d\underline{x}\right)$
is an eigenfunction of
$m_{\cosh_{ \underline{s}}}\circ
\mathcal{K}''_{\underline{c}_{0}}$
with eigenvalue
$\varrho$ iff $Bf$ is an eigenfunction of the operator
$C_{\underline{s}}\circ M_{\underline{\lambda}}$ for the same eigenvalue
$\varrho$.
For
$\underline{s}=\underline{s}_{0}=2 (\pi\beta\underline{J})^
{\frac{1}{\underline{2}}}$,
however, the operator
$C_{\underline{s}_{0}}\circ M_{\underline{\lambda }}:\mathcal{F}_{m}\to\mathcal{F}_{m}$
is conjugate to $\prod^{m}_{i=1}\left( \lambda _{i}\exp \beta J_{i}\right) ^{\frac{1}{2}}\cL_{\beta }
:\cF_{m}^{\underline{\alpha}_{0}^{-1}} \to
 \cF^{\underline{\alpha}_{0}^{-1}}_{m}$ by Proposition \ref{OPCM}.
\qed

Hence one concludes that
$f\in L^{2}\left( \mathbb {R}^{m},dx\right)$
is an eigenfunction of the operator
\[\frac{1}{\prod\limits ^{m}_{l=1}\left( \lambda _{l}\exp \beta J_{l}\right) ^{\frac{1}{2}}}
m_{\cosh_{ \underline{s}_{0}}}\circ \mathcal{K}''_{\underline{c}_{0}}
\]
iff the function $\left(\nu_{\underline{\alpha }_{0}^{ -1}}\circ
B\right)f\in\cF^{\underline{\alpha}_{0}^{-1}}_{m} $ is an
eigenfunction of the operator $\cL_{\beta }$ for the same
eigenvalue. If therefore $f( \underline{x})$ is an eigenfunction of
the operator
\[ \mathcal{\widetilde{G}}_{\beta,\underline{c}_{0}}:=\frac{1}{\prod\limits ^{m}_{l=1}
\left( \lambda _{l}\exp \beta J_{l}\right) ^{\frac{1}{2}}}
  m_{\cosh _{\underline{s}_{0}}}\circ
  \mathcal{K}''_{\underline{c}_{0}}
\]
in
$L^{2}\left( \mathbb {R}^{m},d\underline{x}\right)$,
then the corresponding eigenfunction
$F\left( \underline{z}\right)$
of the operator
$\cL_{\beta }:\cF^{\underline{\alpha }^{-1}_{0}}_{m}\to \cF^{\underline{\alpha }_{0}^{-1}}_{m}$
has the following explicit form
\begin{eqnarray*}
F(z)
& = &\left( \nu_{\underline{\alpha }_{0}^{-1}}\circ B f\right)(\underline{z})\\
& = & 2^{\frac{m}{4}}\int\limits _{\mathbb {R}^{m}}f\left( \underline{x}\right)
     \exp \left( 2\pi \underline{x}\cdot \left( \frac{\beta }{\pi }\underline{J}
     \right) ^{\frac{1}{2}}\circ \underline{z}-\pi \underline{x}^{2}-\frac{\pi }{2}
     \left( \left( \frac{\beta }{\pi }\underline{J}\right) ^{\frac{1}{2}}\circ
     \underline{z}\right) ^{2}\right) d\underline{x}\\
& = & 2^{\frac{m}{4}}\int\limits _{\mathbb {R}^{m}}f\left( \underline{x}\right)
     \exp \left( 2\pi \sum ^{m}_{i=1}\left( x_{i}\sqrt{\frac{\beta J_{i}}{\pi }}z_{i}\right)
     -\pi \underline{x}^{2}-\frac{\beta }{2}\sum ^{m}_{i=1}J_{i}z^{2}_{i}\right)
     d\underline{x}.
\end{eqnarray*}
On the other hand given an eigenfunction
$F=F(z)$ of the operator
$\cL_{\beta }:\cF^{\underline{\alpha }^{-1}_{0}}_{m}\to \cF_{m}^{\underline{\alpha }^{-1}_{0}}$,
the corresponding eigenfunction
$f=f(x)$
of the operator \(\mathcal{\widetilde{G}}_{\beta,\underline{c}_{0}}\)
has the form
\[
f(x)=\left( B^{-1}\circ \nu_{\underline{\alpha}_{0}}F\right)
(x)=2^{\frac{m}{4}} \int\limits _{\mathbb {C}^{m}}F\left(
\underline{\alpha}_{0}\circ \underline{z}\right) \exp \left( 2\pi
\underline{x}\cdot \underline{z}^{*} -\pi
\underline{x}^{2}-\frac{\pi }{2} \underline{z}^{*2}\right)
\exp\left( -\pi \left| \underline{z}\right|
^{2}\right)d\underline{z}.
\]
Inserting the explicit expression for
$\underline{\alpha }_{0}
=\left( \frac{1}{\sqrt{J_{1}}},\ldots ,\frac{1}{\sqrt{J_{m}}}\right)
\sqrt{\frac{\pi }{\beta }}$
we therefore get
\[
f( \underline{x}) =2^{\frac{m}{4}}\int\limits _{\mathbb {C}^{m}}F
\left( \sqrt{\frac{\pi }{\beta J_{1}}}z_{1},\ldots,
\sqrt{\frac{\pi }{\beta J_{m}}}z_{m}\right) \exp \left( 2\pi
\underline{x}\cdot \underline{z}^{*}-\pi
\underline{x}^{2}-\frac{\pi }{2} (\underline{z}^*)^{2}-\pi \left|
\underline{z}\right| ^{2}\right) d\underline{z}
\]
Obviously the integral operator
\(\widetilde{\mathcal{G}}_{\beta,\underline{c}_{0}}\) depends on
\(\beta\) holomorphically and hence defines a holomorphic family
of trace class operators in the Hilbert space
\(L^{2}(\mathbb{R}^{m},d\underline{x} )\). Its Fredholm
determinant \(\det(1-z\mathcal{\widetilde{G}}_{\beta})\) hence is
an entire function in the entire \(\beta\) plane coinciding with
the Fredholm determinant of the Ruelle operator
\(\mathcal{L}_{\beta}\). To relate finally the operator
$m_{\cosh_{ \underline{s}_{0}}}\circ
\mathcal{K}''_{\underline{c}^{0}}$ and its eigenfunctions to those
of the Kac-Gutzwiller operator $\cK_{\beta }$ with kernel
$\cK_{\beta }\left( \xi ,\eta \right)$ we use
\bigskip

\begin{Prop}\label{SpecAequ}
For real $\beta$ consider the two operators
$m_{\cosh_{\underline{s}_{0}}}\circ \mathcal{K}''_{\underline{c}_{0}}$ and
$m_{\sqrt{\cosh_{\underline{s}_{0}}}} \circ \mathcal{K}''_{\underline{c}_{0}}
\circ m_{\sqrt{\cosh_{\underline{s}_{0}}}} $ acting on
the space $L^{2}(\mathbb{R}^{m},d\underline{x}).$
Then the following two statements hold:

\begin{enumerate}
\item[{\rm(i)}]
If $f\in L^{2}(\mathbb{R}^{m},d\underline{x})$ is an eigenfunction of the operator
$m_{\cosh_{\underline{s}_{0}}}\circ \mathcal{K}''_{\underline{c}_{0}}$ with eigenvalue
$\rho \neq 0$, then the function $g := \frac{f}{\sqrt{\cosh_{\underline{s}_{0}}}}$ is an
eigenfunction of the operator
$m_{\sqrt{\cosh_{\underline{s}_{0}}}} \circ \mathcal{K}''_{\underline{c}_{0}}\circ
m_{\sqrt{\cosh_{\underline{s}_{0}}}} $
in $L^{2}(\mathbb{R}^{m},d\underline{x})$  for the same eigenvalue.

\item[{\rm(ii)}]
Conversely,
if $g\in L^{2}(\mathbb{R}^{m},d\underline{x})$ is an eigenfunction of the operator
$m_{\sqrt{\cosh_{\underline{s}_{0}}}} \circ \mathcal{K}''_{\underline{c}_{0}}\circ
m_{\sqrt{\cosh_{\underline{s}_{0}}}} $ with eigenvalue $\rho \neq 0$,
then the function $f=\sqrt{\cosh_{\underline{s}_{0}}}\cdot g $
is an eigenfunction of the operator $m_{\cosh_{\underline{s}_{0}}}\circ
\mathcal{K}''_{\underline{c}_{0}}$ in
$L^{2}(\mathbb{R}^{m},d\underline{x})$ for the same eigenvalue.
\end{enumerate}
\end{Prop}

\pf
\begin{enumerate}
\item[(i)]
If  $f\in L^{2}(\mathbb{R}^{m},d\underline{x})$ is an eigenfunction of
$m_{\cosh_{\underline{s}_{0}}}\circ \mathcal{K}''_{\underline{c}_{0}}$,  then
$g:=\frac{f}{\sqrt{\cosh_{\underline{s}}}}$ is in $L^2(\R^m,d\underline{x})$
since $\cosh_{\underline{s}_0}\ge 1$ and a simple calculation shows that
this function is an eigenfunction of the operator
$m_{\sqrt{\cosh_{\underline{s}_{0}}}} \circ \mathcal{K}''_{\underline{c}_{0}}
\circ m_{\sqrt{\cosh_{\underline{s}_{0}}}}$ for the same eigenvalue.

\item[(ii)]
For $\phi\in L^2(\R^m,\cosh_{\underline{s}_0}(\underline{x})^{-1}d\underline{x})$ set
$$h(\underline{x}):=\int_{\R^m}\tilde{\mathcal{K}}
(2\sqrt\pi \underline{x},2\sqrt\pi \underline{y})\phi(\underline{y})\, d\underline{y}.$$
Then (\ref{KacGutzKernmod}) shows that there exist constants $c,d>0$ such that
$|h(\underline{x})|\le c e^{-d|\underline{x}|^2}$. In particular, it follows that
$\cosh_{\underline{s}_0} h\in L^2(\R^m,d\underline{x})$.

If now $g \in L^{2}(\mathbb{R}^{m},d\underline{x})$ is an eigenfunction
of the operator
$m_{\sqrt{\cosh_{\underline{s}_{0}}}} \circ \mathcal{K}''_{\underline{c}_{0}}
\circ m_{\sqrt{\cosh_{\underline{s}_{0}}}}$ for the nonzero eigenvalue $\rho$,
then a simple calculation shows that the function
$f:=g\sqrt{\cosh_{\underline{s}}}  \in
L^{2}(\mathbb{R}^{m},\frac{1}{\cosh_{\underline{s}}}d\underline{x})$
is an eigenfunction of the operator
$m_{\cosh_{\underline{s}_{0}}}\circ \mathcal{K}''_{\underline{c}_{0}}$
for the eigenvalue $\rho\not=0$ and the above argument applied to
$\phi=\rho (4\pi)^{-\frac{m}{2}}f$
shows that
$f\in L^{2}(\mathbb{R}^{m},d\underline{x})$.
\end{enumerate}
\qed

This leads us to the  main result of this section:

\begin{Thm}\label{THM}
 The Ruelle operator
$\cL_{\beta }:\cF^{\underline{\alpha }^{-1}_{0}}_{m}\to \cF^{\underline{\alpha }^{-1}_{0}}_{m}
$
and the modified Kac-Gutzwiller operator
$\widetilde{\cG}_{\beta }:L^{2}\left( \mathbb {R}^{m},d\underline{\xi} \right) \to
L^{2}\left( \mathbb {R}^{m},d\underline{\xi} \right)$
with kernel
\[\widetilde{\cG}_{\beta }\left( \underline{\xi },\underline{\eta }\right)
=\prod\limits ^{m}_{l=1}\left( \lambda _{l}\exp \beta J_{l}\right) ^{-\frac{1}{2}}
\cosh\left(\sum_{l=1}^{m}\sqrt{\beta J_{l}}\xi_{l}\right)
\widetilde{\mathcal{K}} \left( \underline{\xi },\underline{\eta }\right)\]
with $\widetilde{\mathcal{K}}(\underline{\xi},\underline{\eta})$ defined in (\ref{KacGutzKernmod}) have the same spectrum.
For real \(\beta\) this spectrum coincides also with the spectrum of the
Kac-Gutzwiller operator \(\mathcal{G}_{\beta}\) defined in (\ref{KacGutzmod})
on the Hilbert space
\(L^{2}(\mathbb{R},d\underline{\xi})\).
For nonvanishing eigenvalues the eigenfunctions $F(\underline{z})$ and
$f(\underline{\xi})$ of \(\mathcal{L}_{\beta}\) and \(\mathcal{G}_{\beta}\)
can be related to each other as follows:
\begin{eqnarray*}
f\left( \underline{\xi }\right)
& =&\left( \frac{1}{2\pi }\right) ^{\frac{m}{4}}
\frac{1}{\sqrt{\cosh \sum ^{m}_{i=1}\sqrt{\beta J_{i}}\xi _{i}}}
    \int _{\mathbb {C}^{m}}F\left( \left(
    \frac{\pi }{\beta {J}_1}\right)^{\frac{1}{2}}{z}_1,\ldots,
    \left(\frac{\pi }{\beta {J}_m}\right)^{\frac{1}{2}}{z}_m\right)\cdot\\
&&\cdot    \exp \left( \sqrt{\pi }\underline{\xi }\cdot
\underline{z}^{*}-\frac{1}{4}
    \underline{\xi }^{2}-\frac{\pi }{2}\underline{z}^{*2}
    -\pi \left| \underline{z}\right| ^{2}\right) d\underline{z}\\
F\left( \underline{z}\right)
& =&\left( 8\pi \right) ^{\frac{m}{4}}\int _{\mathbb {R}^{m}}\sqrt{\cosh
\left( 2\sqrt{\pi\beta }
 \sum ^{m}_{l=1}\sqrt{J_{l}}x_{l}\right)}\,
 f\left( 2\sqrt{\pi }\underline{x}\right)
    \cdot\\
&&\cdot    \exp \left( 2\sqrt{\pi \beta }\sum
^{m}_{l=1}\sqrt{J_{l}x_{l}}z_{l}
    -\pi \underline{x}^{2}-\frac{\beta }{2}\cdot \sum ^{m}_{l=1}J_{l}z^{2}_{l}
    \right)
      d\underline{x}.
\end{eqnarray*}
\end{Thm}

 \pf
 We have seen already in Proposition \ref{Lconjugate} that the operators
$\prod ^{m}_{l=1}\left( \lambda _{l}\exp   \beta J_{l}\right) ^{-\frac{1}{2}}
   m_{\cosh_{\underline{s}_{0}}}\circ
   \mathcal{K}''_{\underline{c}_{0}}$
and  $\cL_{\beta }$ are conjugate via the map
$\nu_{\underline{\alpha} ^{-1}_{0}}\circ B$. In fact, we have the commutative diagram
$$
\xymatrix@M=17pt{
L^{2}\left( \mathbb {R}^{m},d\underline{\xi }\right)
 \ar[r]^{\mathcal{K}_{\beta}}
 \ar[d]_{R_{\underline{c}_0}}
& L^{2}\left( \mathbb {R}^{m}d\underline{\xi }\right)
 \ar[d]^{R_{\underline{c}_0}}
&
\\
L^{2}\left( \mathbb {R}^{m},a\left( \underline{x}\right) d\underline{x}\right)
 \ar[d]_{m_{\sqrt{\cosh_{\underline{s}_0}}}}
 \ar[r]
&L^{2}\left( \mathbb {R}^{m},d\underline{x}\right)
 \ar[d]^{m\frac{1}{\sqrt{\cosh_{\underline{s}_0}}}}
&
\\
L^{2}\left( \mathbb {R}^{m},d\underline{x}\right)
 \ar[r]_{\mathcal{K}''_{\underline{c}_0}}
 \ar[d]_{B}
& L^{2}\left( \mathbb {R}^{m},d\underline{x}\right)
  \ar[d]^{B}
  \ar[r]_{m_{\cosh_{\underline{s}_0}}}
& L^{2}\left( \mathbb {R}^{m},d\underline{x}\right)
  \ar[d]^{B}
\\
\mathcal{F}_m
  \ar[r]_{M_{\underline{c}_0}}
  \ar[d]_{\nu_{\underline{\alpha}}}
& \mathcal{F}_m
  \ar[r]_{C_{\underline{s}_0}}
& \mathcal{F}_m
  \ar[d]^{\nu_{\underline{\alpha}}}
\\
\mathcal{F}_m^{\underline{\alpha}^{-1}}
 \ar[rr]_{ \mathcal{L}_\beta \sqrt{\prod_{l=1}^m\lambda_le^{\beta J_l}}}
&
& \mathcal{F}_m^{\underline{\alpha}^{-1}}
}
$$
\bigskip
with
$\underline{c}_0=2\sqrt{\pi}(1,\ldots,1)$ and $\underline{s}_0=2\sqrt{\beta\pi\underline{J}}$.
%
Furthermore
Proposition \ref{SpecAequ} shows that for real $\beta$ the
operators $m_{\cosh_{ \underline{s}_{0}}}\circ
\mathcal{K}''_{\underline{c}_{0}}$ and $m_{\sqrt{\cosh_{
\underline{s}_{0}}}}\circ \mathcal{K}''_{\underline{c}_{0}} \circ
m_{\sqrt{\cosh_{ \underline{s}_{0}}}}$ have the same nonvanishing
eigenvalues. But the last operator is conjugate to the operator
$\mathcal{K}_{\beta }$ with kernel
$\mathcal{K}_{\beta}(\underline{\xi },\underline{\eta })$ through
the map $R_{\underline{c}_{0}}$. Hence if $f\in L^{2}\left(
\mathbb {R}^{m},d\underline{\xi }\right)$ is an eigenfunction of
the integral operator $\cG_{\beta }$, then $\left(
R_{\underline{c}_{0}}f\right) \left( \underline{x}\right)
=\underline{c}_{0}^{\underline{\frac{1}{2}}}f\left(
\boldsymbol{C}_{0}\underline{x}\right)$ is an eigenfunction of the
integral operator $\prod_{l=1}^{m}\left(\lambda_{l}e^{\beta
J_{l}}\right)^{-\frac{1}{2}} m_{\sqrt{\cosh_{
\underline{s}_{0}}}}\circ \mathcal{K}''_{\underline{c}_{0}} \circ
m_{\sqrt{\cosh_{ \underline{s}_{0}}}}$ for the same eigenvalue.
But then
$$
\sqrt{\cosh_{ \underline{s}_{0}}(\underline{x})} \left(
R_{\underline{c}_{0}}f\right) \left( \underline{x}\right) =
\sqrt{\cosh_{ \underline{s}_{0}}(\underline{x})}\quad
\underline{c}_{0}^{\frac{1}{\underline{2}}}f(\boldsymbol{C}_{0}\underline{x})
$$
is an eigenfunction of the operator $\frac{1}{\prod \limits
^{m}_{l=1}\left( \lambda _{l} e^{\beta J_{l}}\right)
^{\frac{1}{2}}}\cdot m_{\cosh_{ \underline{s}_{0}}} \circ
\mathcal{K}''_{\underline{c}_{0}}$ again for the same eigenvalue.
Therefore
\[
F\left( \underline{z}\right) =2^{\frac{m}{4}}\int \limits
_{\mathbb {R}^{m}} \sqrt{\cosh \underline{s}_{0}\cdot
\underline{x}}\quad \underline{c}_{0}^{\frac{1}{\underline{2}}}
f\left( \boldsymbol{C}_{0}\underline{x}\right) \cdot \exp \left(
2\pi \sum ^{m}_{i=1} \sqrt{\frac{\beta J_{i} }{\pi }} x_{i}z_{i}
-\pi \underline{x}^{2}-\frac{\beta }{2}\sum \limits
^{m}_{i=1}J_{i}z^{2}_{i}\right) d\underline{x}
\]
is an eigenfunction of the operator $\cL_{\beta
}:F^{\underline{\alpha }_{0}^{-1}}_{m}\to F^{\underline{\alpha
}^{-1}_{0}}_{m}$ for yet again the  same eigenvalue.

Inserting $\underline{c}_{0}=\left( 2\sqrt{\pi },\ldots
,2\sqrt{\pi }\right)$ one  therefore finds for $F\left( z\right)$
\begin{eqnarray*}
F(\underline{z})
& =&(8\pi)^{\frac{m}{4}}\int _{\mathbb {R}^{m}}
    \sqrt{\cosh \left(2\sqrt{\beta \pi} \sum \limits ^{m}_{l=1}\sqrt{J_{l}}\cdot x_{l}\right)}
    f\left( 2\sqrt{\pi }\underline{x}\right)\cdot \\
&&  \cdot\exp \left( 2\sqrt{\pi \beta }\sum \limits
^{m}_{l=1}\sqrt{J_{l}}x_{l}z_{l}
    -\pi \underline{x}^{2}-\frac{\beta }{2}\sum \limits ^{m}_{l=1}J_{l}z^{2}_{l}
    \right)
    d\underline{x}
\end{eqnarray*}
On the other hand, given an eigenfunction $F\in
\cF^{\underline{\alpha} ^{-1}_{0}}_{m}$ of the operator
$\cL_{\beta }$ we know that $h(\underline{x})=\left( B^{-1}\circ
\nu_{\underline{\alpha }_{0}}\circ F\right) (\underline{x})$ is an
eigenfunction of the operator $\frac{1}{\prod
^{m}_{l=1}\sqrt{\lambda _{l}\exp \beta J_{l}}} m_{\cosh_{
\underline{s}_{0}}}\circ \mathcal{K}''_{\underline{c}_{0}}$.

Then by Proposition \ref{SpecAequ} for real $\beta$ the function
 \(\frac{1}{ \sqrt{\cosh_{ \underline{s}_{0}}(\underline{x})}}h(x) \)
is an eigenfunction of the operator
\[ \frac{1}{\prod_{l=1}^{m}
\left( \lambda _{l}\exp \beta J_{l}\right) ^{\frac{1}{2}}}
m_{\sqrt{\cosh_{ \underline{s}_{0}}}}\circ
\mathcal{K}''_{\underline{c}_{0}}\circ m_{\sqrt{\cosh_{
\underline{s}_{0}}}},
\]
which is again conjugate to \(
\mathcal{G}_{\beta } \) via \( R^{-1}_{\underline{c}_{0}} \) and
hence \( \left( R_{\underline{c}^{-1}_{0}}\left(
\frac{1}{\sqrt{\cosh_{ \underline{s}_{0}}}}\cdot h\right)\right)
\left( \underline{\xi }\right)   \) is an eigenfunction of \(
\cG_{\beta } \). Inserting all the transforms involved we finally
get for the corresponding eigenfunction \( f=f\left(\underline{
\xi} \right) \in L^{2}\left( \mathbb {R}^{m},d\underline{\xi}\right)  \)

\begin{eqnarray*}
f\left( \underline{\xi }\right)
&=&\left( \frac{1}{2\pi }\right) ^{\frac{m}{4}}
\frac{1}{\sqrt{\cosh \left( \sum \limits ^{m}_{l=1}\sqrt{\beta J_{l}}
\xi _{l}\right) }}\cdot \int \limits _{\mathbb {C}^{m}}F\left(
\sqrt{\frac{\pi }{\beta J_{1}}}z_{1},\ldots ,\sqrt{\frac{\pi }{\beta J_{m}}}z_{m}\right) \cdot\\
&&\cdot\exp \left( \sqrt{\pi }\underline{\xi }\cdot \underline{z}^{*}-\frac{1}{4}\underline{\xi }^{2}
-\frac{\pi }{2}\underline{z}^{*2}-\pi \left| \underline{z}\right| ^{2}\right) d\underline{z}
\end{eqnarray*}
\qed

To discuss the zeros and poles of the Ruelle zeta function for the Kac-Baker models we 
need the following proposition.

\begin{Prop}\label{REAL}
For $\beta$ real the operators $\widetilde{\mathcal{G}}_{\beta}$ and $\cL_{\beta}$ have real spectrum.
\end{Prop}
\pf
For $\beta$ real and $f\in L^{2}\left(\mathbb{R}^{m},d\underline{\xi}\right)$ an eigenfunction of 
the operator $\widetilde{\mathcal{G}}_{\beta}$ with eigenvalue $\varrho$ 
consider the scalar product  
\[\left(\frac{f}{\cosh_{\underline{R}_{0}}},\widetilde{\mathcal{G}}_{\beta} f\right)
= \oline{\varrho} \left(\frac{f}{\cosh_{\underline{R}_{0}}},f\right),\] where $\underline{R}_{0}
=\sqrt{\beta \underline{J}}$.
 But if
\[\int_{\mathbb{R}^{m}}\prod_{l=1}^{m}\left(\lambda_{l}\exp{\beta J_{l}}\right)^{-\frac{1}{2}} 
\cosh_{\underline{R}_{0}}(\underline{\xi})\, \widetilde{\mathcal{K}}
\left(\underline{\xi},\underline{\eta}\right)\,f(\underline{\eta})\,d\underline{\eta}
=\varrho f(\underline{\xi}),\]
then
\[\int_{\mathbb{R}^{m}}\prod_{l=1}^{m}\left(\lambda_{l}\exp{\beta J_{l}}\right)^{-\frac{1}{2}} 
\cosh_{\underline{R}_{0}}(\underline{\xi})\, \widetilde{\mathcal{K}}
\left(\underline{\xi},\underline{\eta}\right)\,\frac{f(\underline{\xi})}{\cosh_{\underline{R}_{0}}(\underline{\xi})}
\,d\underline{\xi}
=\varrho\, \frac{f(\underline{\eta})}{\cosh_{\underline{R}_{0}}(\underline{\eta})}\] 
since $\widetilde{\mathcal{K}}\left(\underline{\xi},\underline{\eta}\right)=
\widetilde{\mathcal{K}}\left(\underline{\eta},\underline{\xi}\right)$ 
and the function $\frac{f(\underline{\xi})}{\cosh_{\underline{R}_{0}}(\underline{\xi})}$ 
belongs to the space $L^{2}\left(\mathbb{R}^{m},d\underline{\xi}\right)$ if $f$ is an eigenfunction 
of $\widetilde{\mathcal{G}}_{\beta}$. 
Since $\widetilde \cG$ is real valued we can now calculate
\begin{eqnarray*}
\left(\frac{f}{\cosh_{\underline{R}_{0}}},\widetilde{\mathcal{G}}_{\beta} f\right)
&=&\int_{\mathbb{R}^{m}}\frac{f(\underline{\xi})}{\cosh_{\underline{R}_{0}}(\underline{\xi})}
   \overline{
   \int_{\mathbb{R}^{m}}\widetilde{\mathcal{G}}_{\beta}(\underline{\xi},\underline{\eta})
    f(\underline{\eta})\, d\underline{\eta}
            }
     \, d\underline{\xi}\\
&=&\int_{\mathbb{R}^{m}}\int_{\mathbb{R}^{m}}\frac{f(\underline{\xi})}{\cosh_{\underline{R}_{0}}(\underline{\xi})}
                         \widetilde{\mathcal{G}}_{\beta}(\underline{\xi},\underline{\eta})
   \, d\underline{\xi}
   \overline{f(\underline{\eta})}  \, d\underline{\eta}\\
&=&\left(\varrho\,\frac{f}{\cosh_{\underline{R}_{0}}},f\right)\\
&=&\varrho\left(\frac{f}{\cosh_{\underline{R}_{0}}},f\right)
\end{eqnarray*}
On the other hand
$\left(\frac{f}{\cosh_{\underline{R}_{0}}},f\right)\neq 0$, 
since according to Mehlers formula in Proposition \ref{KacKernProp}
\[\int_{\mathbb{R}^{m}}\int_{\mathbb{R}^{m}}\frac{f(\underline{\xi})}{\cosh_{\underline{R}_{0}}(\underline{\xi})} 
\widetilde{\mathcal{G}}_{\beta}(\underline{\xi},\underline{\eta})\,\overline{f(\underline{\eta})}
\,d\underline{\xi}\,d\underline{\eta}
=2\prod_{l=1}^{m}\left(4 \pi \exp{\beta J_{l}} \right)^{-\frac{1}{2}}
\sum_{\underline{\alpha}\in \mathbb{N}_{0}^{m}}
\underline{\lambda}^{\underline{\alpha}}\left|
\int_{\mathbb{R}^{m}}h_{\underline{\alpha}}
\left(\frac{\underline{\xi}}{2\sqrt{\pi}}\right)\,
\overline{f(\underline{\xi})}
\,d\underline{\xi}\right|^{2}>\,0 
\] 
for $f$ an eigenfunction of $ \widetilde{\mathcal{G}}_{\beta}$. But then it is clear that $\varrho$ must 
be identical to $\overline{\varrho}$ and hence real.
\qed

\section{The Ruelle zeta function for the Kac-Baker model}
\label{ZRRF}

Recall from Proposition \ref{zetaRmero} that the Ruelle zeta
function $\zeta _{R}\left( z,\beta \right) =\exp \sum ^{\infty
}_{n=1}\frac{z^{n}}{n}Z_{n}(\beta )$ with $Z_{n}(\beta )$ the
partition function for the Kac-Baker model can be written as \[
\zeta _{R}\left( z,\beta \right) =\prod _{\underline{\alpha }\in
\left\{ 0,1\right\} ^{m}} \det \left( 1-z\underline{\lambda
}^{\underline{\alpha }}\mathcal{L}_{\beta }\right) ^{\left(
-1\right) ^{\left| \underline{\alpha } \right| +1}} \] with
the operator
$\mathcal{L}_{\beta }:\mathcal{B}(D)\to\mathcal{B}(D)$
given by

 \[ \left( \mathcal{L}_{\beta }f\right)
\left( \underline{z}\right) =e^{\beta \underline{J}\cdot
\underline{z}} f\left( \underline{\lambda }+\boldsymbol{\Lambda}
\underline{z}\right) + e^{-\beta \underline{J}\cdot
\underline{z}}f\left( -\underline{\lambda }+\boldsymbol{\Lambda}
\underline{z}\right), \] 

 where $\boldsymbol{\Lambda}$ is the
$m\times m$-diagonal matrix with the $\lambda_j$, $j=1,\ldots,m$,
as diagonal elements. 

\begin{Rem}\label{specL0}
For $\beta=0$ one then finds for the
spectrum $\sigma\left(\mathcal{L}_{0}\right)$ of $\cL_0$
$$\sigma \left(
\cL_{0}\right) =
 \left\{ 2\underline{\lambda }^{\underline{\alpha }}\mid
  \underline{\alpha }=\left( \alpha _{1},\ldots ,\alpha _{m}\right)
  \in \mathbb {N}^{m}_{0}\right\}.
  $$
In fact, by induction over $|\underline{\alpha}|$ one finds a polynomial
eigenfunction for each $2\underline{\lambda }^{\underline{\alpha }}$.
In the case $m=1$ one obtains $1, z, z^2+\frac{\lambda^2}{\lambda^2-1}, 
z^3+ \frac{3\lambda^2}{\lambda^2-1}z$ and so on. On the other hand, one can show that
all eigenfunctions are polynomial, since the derivative of an eigenfunction for the 
eigenvalue $\rho$ is again an eigenfunction for the eigenvalue $\rho\lambda^{-1}$
(in the case $m=1$) so that infinitely non-vanishing derivatives of an eigenfunction
would contradict the compactness of the operator. For $m>1$ the argument is similar.
\qed
\end{Rem}

Hence the Ruelle function $\zeta_R \left( z,\beta \right)$  for
$\beta =0$ has the form
\[
\zeta_R \left( z,0\right) =\prod _{\underline{\alpha }\in \left\{
0,1\right\} ^{m}} \prod _{\underline{\beta }\in \mathbb
{N}^{m}_{0}}\left( 1-2z\underline{\lambda }^{\underline{\alpha
}+\underline{\beta }} \right) ^{\left( -1\right) ^{\left|
\underline{\alpha }\right| +1}}.
\]
By induction on $m$ one shows
\[
\zeta_R \left( z,0\right) =\frac{1}{1-2z}.\] Indeed for $m=1$ one
has
\[
\zeta_R \left( z,0\right) =\frac{\det \left( 1-z\lambda
\mathcal{L}_{0}\right) }{\det \left( 1-z\mathcal{L}_{0}\right) }
=\frac{\prod ^{\infty }_{n=0}\left( 1-2z\lambda ^{n+1}\right)
}{\prod ^{\infty }_{n=0}\left( 1-2z\lambda ^{n}\right) }
=\frac{1}{1-2z}.
\]
For $m=n+1$,   on the other hand, one calculates
\begin{eqnarray*}
\zeta_R \left( z,0\right)
& = & \prod _{\alpha _{n+1}\in \left\{ 0,1\right\} }
      \prod _{\underline{\alpha }\in \left\{ 0,1\right\} ^{n}}
      \prod _{\beta _{n+1}\in \mathbb {N}_{0}}
      \prod _{\underline{\beta }\in \mathbb {N}_{0}^{n}}
      \left( 1-2z\lambda ^{\alpha _{n+1}+\beta _{n+1}}_{n+1}
      \underline{\lambda }^{\underline{\alpha }
      +\underline{\beta }}\right) ^{\left( -1\right) ^{1+\left|
      \underline{\alpha }\right| +\alpha _{n+1}}}\\
& = & \prod _{\underline{\alpha }\in \left\{ 0,1\right\} ^{n}}
      \prod _{\beta _{n+1}\in \mathbb {N}_{0}}
      \prod _{\underline{\beta }\in \mathbb {N}^{n}_{0}}
      \left( 1-2z\lambda ^{\beta _{n+1}}_{n+1}
      \underline{\lambda }^{\underline{\alpha }
      +\underline{\beta }}\right) ^{\left( -1\right) ^{1+\left|
      \underline{\alpha }\right| }}\cdot\\
&& \cdot \prod _{\underline{\alpha }\in \left\{ 0,1\right\} ^{n}}
         \prod _{\beta _{n+1}\in \mathbb {N}_{0}}
         \prod _{\underline{\beta }\in \mathbb {N}^{n}_{0}}
         \left( \left( 1-2z\lambda ^{\beta _{n+1}+1}_{n+1}
         \underline{\lambda }^{\underline{\alpha }
         +\underline{\beta }}\right) ^{\left( -1\right) ^{1+\left|
         \underline{\alpha }\right| -1}}\right) \\
& = & \prod _{\beta _{n+1}\in \mathbb {N}_{0}}
      \left( \frac{1}{1-2z\lambda ^{\beta _{n+1}}_{n+1}}\right)
      \prod _{\beta _{n+1}\in \mathbb {N}_{0}}
      \left( 1-2z\lambda ^{\beta _{n+1}+1}_{n+1}\right) \\
& = & \frac{1}{1-2z} .
\end{eqnarray*}

This result does not come as a surprise since for   $\beta =0$ the
Ruelle zeta function for the Kac-Baker model is just the
Artin-Mazur zeta function for the full subshift over two symbols
determined by Bowen and Lanford in \cite{BoLa75}.

 To determine the zeros and poles of the Ruelle function $\zeta_{R}(\beta):=\zeta_{R} \left( 1,\beta\right)$ 
on the real $\beta$-axis one has to investigate the zeros
of the Fredholm determinants $\det \left(
1-\underline{\lambda}^{\underline{\alpha }}\mathcal{L}_{\beta
}\right)$. Obviously this function takes the special value $-1$ at the point $\beta = 0$. We know already from our discussion of the
Kac-Gutzwiller operator $\mathcal{G}_{\beta }$ that all its eigenvalues and
hence also those of the Ruelle operator
$\mathcal{L}_{\beta }$ are real for real $\beta$ and nonnegative for $\beta \geqslant 0$. The poles and zeros of the
Ruelle zeta function $\zeta _{R}\left( 1,\beta \right)$ can be
located only at those values of $\beta$  where one of the numbers $\left(
\underline{\lambda }^{\underline{\alpha }}\right)
^{-1},\underline{\alpha }\in \left\{ 0,1\right\} ^{m}$ belongs to
the spectrum $\sigma \left( \mathcal{L}_{\beta }\right)$ of
$\mathcal{L}_{\beta }$. For   $\beta =0$   we have seen that
$\sigma \left( \mathcal{L}_{0}\right) =2\underline{\lambda
}^{\underline{\beta }}$ with $\underline{\beta }\in \mathbb
{N}_{0}^{m}$. But
$$
\left( \underline{\lambda }^{\underline{\alpha }}\right) ^{-1}
=2\underline{\lambda }^{\underline{\beta }}
\quad\Leftrightarrow\quad
 \frac{1}{2} =\underline{\lambda}^{\underline{\alpha }+\underline{\beta }}
$$
can be true only for finitely many $\underline{\beta }\in \mathbb
{N}^{m}_{0}$
    since
 $\underline{\alpha }\in \left\{ 0,1\right\} ^{m} $
    can take only finitely many values.

Furthermore, for infinitely many
$\underline{\beta }\in \mathbb {N}^{m}_{0}$
    the eigenvalue
 $2\underline{\lambda }^{\underline{\beta }}$
 is strictly smaller than
$\left( \underline{\lambda }^{\underline{\alpha }}\right) ^{-1}$
    for all
 $\underline{\alpha }\in \left\{ 0,1\right\} ^{m}$.
 Hence, if we can show that infinitely many eigenvalues
 $\varrho (\beta )$
 of
$\cL_{\beta }$
    tend to
$+\infty$
    for
$\left| \beta \right| \to \infty$,
   then the Ruelle zeta function will have infinitely many
  ``nontrivial''
 zeros and poles on the real
$\beta$-axis at least for generic values of $\underline{\lambda}$ for which possible
cancellations of zeros in the quotients of the different Fredholm determinants do not occur.

To derive the asymptotic behavior of the eigenvalues of the
transfer operator for $\left| \beta \right| \to \infty$ we remark
first that the eigenspace of any eigenvalue $\varrho$ of
$\mathcal{L}_{\beta}$ has a basis consisting of eigenfunctions
which are also eigenfunctions of the operator
$P:\mathcal{B}(D)\to\mathcal{B}(D)$ defined as
\[
Pf\left( \underline{z}\right) =f\left( -\underline{z}\right). \]
 We call the eigenfunctions with
$Pf=f$ even and those with $Pf=-f$ odd. Consider then the two
operators $\mathcal{L}^{\pm }_{\beta }:\mathcal{B}(D)\to
\mathcal{B}(D)$ defined via
\[
\cL^{+}_{\beta }f\left( \underline{z}\right) =e^{\beta
\underline{J}\cdot \underline{z}} f\left( \underline{\lambda
}+\boldsymbol {\Lambda }\underline{z}\right) +e^{-\beta
\underline{J}\cdot \underline{z}} f\left( \underline{\lambda
}-\boldsymbol {\Lambda }\underline{z}\right),
\]
respectively,
\[
\mathcal{L}^{-}_{\beta }f\left( \underline{z}\right) =e^{\beta
\underline{J}\underline{z}}f\left( \underline{\lambda
}+\boldsymbol {\Lambda }\underline{z}\right) -e^{-\beta
\underline{J}\cdot \underline{z}}f\left( \underline{\lambda }
-\boldsymbol {\Lambda }\underline{z}\right). \]
 The
eigenfunctions of $\mathcal{L}^{+}_{\beta }$ and
$\mathcal{L}^{-}_{\beta }$ are just the even, respectively odd,
eigenfunctions of $\mathcal{L}_{\beta }$. We call the
corresponding eigenvalues even, respectively odd.
Since $\mathcal{L}_\beta$ and $P$ commute
one therefore finds for the spectrum $\sigma( \mathcal{L}_{\beta
})$ of $\mathcal{L}_{\beta }$:
\[
\sigma ( \mathcal{L}_{\beta }) =\sigma ( \mathcal{L}^{+}_{\beta })
\cup \sigma ( \mathcal{L}^{-}_{\beta })
\]
Extending a result of B.~Moritz (see \cite{Mo89}) for $m=1$ to the general case
$m\in \mathbb {N}$ we find for the restricted range $0 \leqslant \lambda_{i} \leqslant \frac{1}{2}$ for $1\leqslant i \leqslant m$ of the parameters $\underline{\lambda}$ and arbitrary $N \geqslant 1$:

\begin{Thm}\label{zetaasympThm}
\begin{enumerate}
\item[{\rm(i)}] For $\beta\to\pm\infty$ the $N$ leading even
eigenvalues $\rho_{\underline{\alpha}}$ of $\cL_{\beta }$ behave like
$$ \underline{\lambda }^{\underline{\alpha }}
\left( \pm 1\right)^{\left|\underline{\alpha }\right|} \exp
\left(\beta \underline{J}\cdot\left( 1-\boldsymbol {\Lambda }\right)
^{-1}\underline{\lambda }\right).
$$

\item[{\rm(ii)}] For $\beta \to \pm\infty$ the $N$ leading odd
eigenvalues $\rho_{\underline{\alpha}}$ of $\cL_{\beta}$ behave like
$$\underline{\lambda }^{\underline{\alpha }} \left( \pm 1\right)
^{\left|\underline{\alpha }\right| +1}\exp \left(\beta
\underline{J}\cdot \left( 1-\boldsymbol {\Lambda }\right)
^{-1}\underline{\lambda }\right).
$$
\end{enumerate}
\end{Thm}

\pf
 Consider first the even eigenvalues and their asymptotic behavior for
$\beta \to +\infty$ .
 Then
$$\mathcal{L}^{+}_{\beta }g\left( \underline{z}\right)
=\exp\left( \beta \underline{J}\cdot \underline{z}\right)
  g\left(\underline{\lambda } +\boldsymbol {\Lambda }\underline{z}\right)
+\exp\left( -\beta \underline{J}\cdot \underline{z}\right)
 g\left(\underline{\lambda }-\boldsymbol {\Lambda }\underline{z}\right)
=\varrho\, g(z).$$
Writing $g(z)=\exp\left( \beta
\underline{J}\cdot \left( 1-\boldsymbol {\Lambda }\right)
^{-1}\underline{z}\right) u( \underline{z}) $
 we get for
$u$ the equation
\begin{eqnarray*}
&&\hskip -4em u\left(
\underline{\lambda}+\boldsymbol{\Lambda}\underline{z}\right)
+\exp\left( -2\beta \underline{J}\cdot\left( 1-\boldsymbol
{\Lambda }\right) ^{-1}\boldsymbol{\Lambda}\underline{z}\right)
\exp\left( -2\beta \underline{J}\cdot \underline{z}\right) u\left(
\underline{\lambda }-\boldsymbol {\Lambda }\underline{z}\right)
=\\
&=&\varrho \, \exp\left( -\beta \underline{J}\cdot \left( 1-\boldsymbol
{\Lambda }\right) ^{-1}\underline{\lambda }\right)u(\underline{z})
\end{eqnarray*}
and hence
\[
u\left( \underline{\lambda }+\boldsymbol {\Lambda
}\underline{z}\right) +\exp\left( -2\beta \underline{J}\cdot \left(
1-\boldsymbol {\Lambda }\right) ^{-1}\underline{z}\right)
u(\underline{\lambda }-\boldsymbol {\Lambda }\underline{z})
=\overline{\varrho } u( \underline{z}), \] where
\[
\overline{\varrho }=\varrho \exp\left( -\beta \underline{J}\cdot \left(
1-\boldsymbol {\Lambda }\right) ^{-1}\underline{\lambda
}\right).\]
Replacing $\underline{z}$ by
$\underline{\lambda}+\underline{z}$ and introducing the function
$h( \underline{z}) :=u( \underline{\lambda }+\underline{z})$ one
arrives at the equation
\begin{eqnarray*}
\overline{\varrho} \, h\left( \underline{z}\right)
&=&h( \boldsymbol {\Lambda }\underline{\lambda }
     +\boldsymbol {\Lambda }\underline{z})
     +\exp\left( -2\beta \underline{J}\cdot \left( 1-\boldsymbol {\Lambda }\right) ^{-1}
     \underline{\lambda }\right)
\exp\left( -2\beta \underline{J}\cdot\left( 1-\boldsymbol
    {\Lambda }\right) ^{-1}\underline{z}\right) h( -\boldsymbol
    {\Lambda }\underline{\lambda }-\boldsymbol
    {\Lambda}\underline{z}).
\end{eqnarray*}
Define an operator
$\mathcal{T}^+_{\beta;\underline{\lambda}}:\mathcal{B}\left(
D_{\widetilde{\underline{R}}}\right) \to \mathcal{B}\left(
D_{\widetilde{\underline{R}}}\right)$ for $D_{\widetilde{\underline{R}}}$ the polydisc with
$\widetilde{R}_{i}>\frac{\lambda ^{2}_{i}}{1-\lambda _{i}}$ for $1 \leqslant i \leqslant m$ via
\[
\mathcal{T}^{+}_{\beta ;\underline{\lambda }}h( \underline{z})
:=h( \boldsymbol {\Lambda }\underline{\lambda }+\boldsymbol
{\Lambda }\underline{z}) +\exp\left( -2\beta \underline{J}\cdot \left(
1-\boldsymbol {\Lambda }\right) ^{-1}\underline{\lambda}\right)
\exp\left( -2\beta \underline{J}\cdot \left( 1-\boldsymbol {\Lambda
}\right) ^{-1}\underline{z}\right) h( -\boldsymbol {\Lambda
}\underline{\lambda }-\boldsymbol {\Lambda }\underline{z}).
\]
Then $\mathcal{T}^{+}_{\beta ,\underline{\lambda }}$ is nuclear
and its eigenvalues are just the numbers $\overline{\varrho
}=\varrho \exp\left( -\beta \underline{J}\cdot \left( 1-\boldsymbol
{\Lambda }\right) ^{-1}\underline{\lambda }\right)$, where
$\varrho$ is an eigenvalue of the operator $\cL^{+}_{\beta }$. If
$\mathcal{T}_{\underline{\lambda }}:\mathcal{B}\left(
D_{\widetilde{\underline{R}}}\right) \to \mathcal{B}\left(
D_{\widetilde{\underline{R}}}\right)$ denotes then the composition operator
\[
\mathcal{T}_{\underline{\lambda }}h\left( \underline{z}\right) =
h(\boldsymbol {\Lambda }\underline{\lambda }+\boldsymbol {\Lambda
}\underline{z}),
\]
one finds for $\underline{\lambda }=\left( \lambda
_{1},\ldots ,\lambda _{m}\right)$ with $0<\lambda
_{i}<\frac{1}{2}$ for all
 $1\leq i\leq m$:
\[
\lim _{\beta \to +\infty }\left\Vert \mathcal{T}^+_{\beta
,\underline{\lambda }}-\mathcal{T}_{\underline{\lambda
}}\right\Vert =0\] since for these values of $\lambda _{i}$ one
can find $\widetilde{R}_{i}$ such that for all $1\leqslant i \leqslant m$
\[
\lambda _{i}>\widetilde{R}_{i}>\frac{\lambda ^{2}_{i}}{1-\lambda
_{i}}.
\]
But the eigenvalues of the operator
$\mathcal{T}_{\underline{\lambda }}$ can be determined explicitly:
they are given by the numbers $\underline{\lambda
}^{\underline{\alpha }}$ with $\underline{\alpha }\in \N^{m}_{0}.$
>From this the asymptotic behavior of the leading eigenvalues $\varrho_{\underline{\alpha}}$ of
$\mathcal{L}^{+}_{\beta }$ and hence the asymptotic behavior of the leading
even eigenvalues of $\mathcal{L}_{\beta }$ follows immediately.

The proof of the behavior of the odd eigenvalues for $\beta \to
+\infty$ follows the same line of arguments applied to the operator
$\mathcal{L}^{-}_{\beta }$ instead of $\mathcal{L}^{+}_{\beta }$.

For the asymptotic behavior of the even eigenvalues for
   $\beta \to -\infty$
 consider the operator
\[
\widetilde{\mathcal{L}}^{+}_{\beta }g(\underline{z}) =\exp
\left(-\beta \underline{J}\cdot \underline{z}\right)g(
\underline{\lambda }+\boldsymbol {\Lambda }\underline{z})
+\exp\left( \beta \underline{J}\cdot \underline{z}\right) g(
\underline{\lambda }-\boldsymbol {\Lambda }\underline{z})
 \]
and the behavior of its eigenvalues for $\beta \to +\infty$. In
this case we write an eigenfunction $g$ as
\[
g( \underline{z}) =\exp\left( \beta \underline{J}\cdot \left(
1+\boldsymbol {\Lambda }\right) ^{-1}\underline{z}\right) u(
\underline{z})
 \]
and get for $u$ the equation
\[
\exp\left( -2\beta \underline{J}\cdot \underline{z}\right)
\exp\left( 2\beta \underline{J}\cdot \left( 1+\boldsymbol {\Lambda
}\right) ^{-1}\boldsymbol {\Lambda }\underline{z}\right)
 u(\underline{\lambda }+\boldsymbol {\Lambda }z)
 +u(\underline{\lambda }-\boldsymbol {\Lambda }\underline{z})
=\overline{\varrho }u(\underline{z}),
 \]
where $\overline{\varrho }=\varrho \exp \left( -\beta
\underline{J}\cdot \left( 1+\boldsymbol {\Lambda }\right)
^{-1}\underline{\lambda }\right)$.  Introducing the function
$h(
\underline{z}) :=u( \underline{\lambda }+\boldsymbol
{\Lambda}\underline{z})$
we finally arrive at
\[
\overline{\varrho }\, h(\underline{z}) =h( -\boldsymbol {\Lambda
}\underline{\lambda}-\boldsymbol {\Lambda }\underline{z})
+\exp\left( -2\beta \underline{J}\cdot \left( 1+\boldsymbol {\Lambda
}\right) ^{-1}\underline{\lambda}\right)
 \exp \left(-2\beta
\underline{J}\cdot \left( 1+\boldsymbol {\Lambda }\right)
^{-1}\underline{z}\right)h( \boldsymbol {\Lambda
}\underline{\lambda }+\boldsymbol {\Lambda }\underline{z}).
\]
Hence $\overline{\varrho }$ is an eigenvalue of the operator $
\mathcal{T}^-_{\beta;\underline{\lambda }}:\mathcal{B}\left(
D_{\widetilde{\underline{R}}}\right) \to \mathcal{B}\left(
D_{\widetilde{\underline{R}}}\right)
$
 with
\begin{eqnarray*}
\mathcal{T}^-_{\beta;\underline{\lambda }}h( \underline{z}) &=&h(
-\boldsymbol {\Lambda }\underline{\lambda } -\boldsymbol {\Lambda
}\underline{z}) +\exp\left( -2\beta \underline{J}\cdot \left(
1+\boldsymbol {\Lambda }\right) ^{-1}
\underline{\lambda}\right)
\exp\left( -2\beta \underline{J}\cdot \left( 1+\boldsymbol
{\Lambda }\right) ^{-1}\underline{z}\right) h( \boldsymbol
{\Lambda }\underline{\lambda}+\boldsymbol {\Lambda
}\underline{z}).
\end{eqnarray*}
 For $0<\lambda _{i}<\frac{1}{2}$ , $1\leqslant i \leqslant m$ we can find again $\widetilde{R}_{i}$
with $\lambda _{i}>\widetilde{R}_{i}>\frac{\lambda ^{2}_{i}}{1-\lambda _{i}}$
and hence
\[
\lim _{\beta \to +\infty }\left\Vert
\mathcal{T}^-_{\beta;\underline{\lambda }} -\mathcal{T}_{\lambda
}\right\Vert =0,\] where
\[
\mathcal{T}_{\underline{\lambda }}h\left( \underline{z}\right) =h(
-\boldsymbol {\Lambda }\underline{\lambda}-\boldsymbol {\Lambda
}\underline{z}),
 \]
is nuclear of order zero on the Banach space $\mathcal{B}\left(
D_{\widetilde{\underline{R}}}\right)$.
 The spectrum of
 $\mathcal{T}_{\underline{\lambda }}$,
however, is given  by the numbers $\left( -\underline{\lambda
}\right) ^{\underline{\alpha }}$ with $\underline{\alpha }\in \mathbb
{N}^{m}_{0}$
 and hence the leading even eigenvalues
$\varrho_{\underline{\alpha}}$ of the operator $\mathcal{L}_{\beta }$
 behave for
$\beta \to -\infty$
 like
$\left( -1\right) ^{\left| \underline{\alpha }\right|
}\underline{\lambda }^{\underline{\alpha }}\exp
\left(-\beta \underline{J}\cdot \left( 1+\boldsymbol{\Lambda} \right)
^{-1}\underline{\lambda }\right)$.

In exactly the same way one shows that the leading odd eigenvalues
$\varrho_{\underline{\alpha}}$ of $\mathcal{L}_{\beta }$ behave for $\beta \to -\infty$
like $-\left( -1\right) ^{\left| \underline{\alpha }\right|
}\underline{\lambda }^{\underline{\alpha }} \exp\left( -\beta
\underline{J}\cdot\left( 1+\boldsymbol {\Lambda }\right)
^{-1}\underline{\lambda }\right)$. \qed

Theorem \ref{zetaasympThm} shows that both for $\beta \to \pm
\infty$ infinitely many eigenvalues of $\mathcal{L}_{\beta }$
 tend to
$ +\infty$.
 Therefore the determinants
$\det \left( 1-\underline{\lambda }^{\underline{\alpha }}
\mathcal{L}_{\beta }\right)$ with $\underline{\alpha }\in \left\{
0,1\right\} ^{m}$ have infinitely many zeros in the real variable
$\beta$.
 Therefore for generic
$\underline{\lambda }$ with $0<\lambda _{i}<\frac{1}{2}$ for
$1 \leqslant i\leqslant m $
 the Ruelle zeta function has infinitely many poles and zeros on the real
 axis.

 For special values of the parameters
$\underline{\lambda}$ some eigenfunctions and their eigenvalues for
the Ruelle operator $\cL_{\beta}$ are explicitly known. If
$\underline{\lambda}=\underline{\lambda}_{0}$ with
$\lambda_{0,i}=\lambda_{0}=\frac{1}{2}$ for all $1\leqslant i \leqslant m$,
then the  identity

$$\sinh\Big(2\beta\underline{J}\cdot
\big(\frac{1}{2}(\underline{z}\pm\underline{1})\big)\Big) =
\frac{1}{2}\left(e^{\beta\underline{J}\cdot\underline{z}\pm
                 \beta\underline{J}\cdot\underline{1}}-
                 e^{-\beta\underline{J}\cdot\underline{z}\mp
                 \beta\underline{J}\cdot\underline{1}}
           \right)$$

and the calculation

\begin{eqnarray*} &&\hskip -4em
e^{\beta\underline{J}\cdot\underline{z}}
\sinh\Big(2\beta\underline{J}\cdot
\big(\frac{1}{2}(\underline{z}+\underline{1})\big)\Big)+
e^{-\beta\underline{J}\cdot\underline{z}}
\sinh\Big(2\beta\underline{J}\cdot
\big(\frac{1}{2}(\underline{z}-\underline{1})\big)\Big)=\\
&=&\frac{1}{2}e^{\beta\underline{J}\cdot\underline{z}}
   \left(e^{\beta\underline{J}\cdot\underline{z}+
                 \beta\underline{J}\cdot\underline{1}}-
                 e^{-\beta\underline{J}\cdot\underline{z}-
                 \beta\underline{J}\cdot\underline{1}}
           \right)+
\frac{1}{2}e^{\beta\underline{J}\cdot\underline{z}}
   \left(e^{\beta\underline{J}\cdot\underline{z}-
                 \beta\underline{J}\cdot\underline{1}}-
                 e^{-\beta\underline{J}\cdot\underline{z}+
                 \beta\underline{J}\cdot\underline{1}}
           \right)\\
&=&\frac{1}{2}
   \left(e^{2\beta\underline{J}\cdot\underline{z}+
                 \beta\underline{J}\cdot\underline{1}}-
                 e^{-\beta\underline{J}\cdot\underline{1}}
           \right)+
\frac{1}{2}
   \left(e^{-\beta\underline{J}\cdot\underline{1}}-
                 e^{-2\beta\underline{J}\cdot\underline{z}+
                 \beta\underline{J}\cdot\underline{1}}
           \right)\\
&=&\frac{1}{2}e^{\beta\underline{J}\cdot\underline{1}}
   \left(e^{2\beta\underline{J}\cdot\underline{z}}-
                 e^{-2\beta\underline{J}\cdot\underline{z}}
           \right)\\
&=&e^{\beta\underline{J}\cdot\underline{1}}
   \sinh(2\beta\underline{J}\cdot\underline{z})
\end{eqnarray*}

 shows that  the functions
 \begin{equation}\label{EFn}
f_{1,n}(\underline{z})=
P_{n}(\underline{z}) \sinh (2\beta
\underline{J}\cdot\underline{z}),
\end{equation}
with $P_{n}(\underline{z})$ a
polynomial homogeneous of degree $n$ in all the variables $z_{l}$
which is invariant under all translations of the form
 $z_{l} \to z_{l}+ c$ for real $c$, are indeed eigenfunctions of
 $\cL_{\beta}$ with eigenvalue
 $$\varrho_{1,n}=e^{\beta \sum_{l=1}^{m}J_{l}}\lambda_{0}^{n}:=\varrho_1\lambda_0^n.$$
The dimension of the space $\cB(D)_n$ of eigenfunctions defined by (\ref{EFn})
for fixed $n$ is the dimension of the 
above space of polynomials and will be calculated in the following proposition.

\begin{Prop} \label{polydim}
Let $\mathbb{K}$ be $\R$ or $\C$. The dimension of
the space $V_{m,n,\underline{c}}$ of homogeneous polynomials
$f(z_1,\ldots,z_m)\in \mathbb{K}[z_1,\ldots,z_m]$ of degree $n$ in
$m$ variables which are invariant under a fixed non-zero
translation $\underline{z}\mapsto\underline{z}+\underline{c}$ with
$\underline{c}\in \mathbb{K}^n$ is
$$\begin{cases}\textstyle\begin{pmatrix}m+n-2\\
n\end{pmatrix}&\mbox{ for } \ m\ge 2\\
1& \mbox{ for }\  m=1 \mbox{ and } n=0\\
0&\mbox{ otherwise}. \end{cases} $$
 \end{Prop}

\pf Let $f\in V_{m,n,\underline{c}}$. For fixed $\underline z$
consider first the polynomial $t\mapsto
f(\underline{z}-t\underline{c})-f(\underline{z})$. It is of degree
less or equal $n$ and has infinitely many zeros (for $t\in \Z$),
so it is zero. Thus we have
$$f(\underline{z}+\mathbb{K}\underline{c})=f(\underline{z})\quad
\forall \underline{z}\in \mathbb{K}^m.$$ Changing coordinates we
may now assume that $\underline{c}=(1,0,\ldots,0)$, i.e. $f$
depends only on the variables $z_2,\ldots,z_m$.  Thus for $m=1$
only the constant polynomials satisfy the required invariance,
whereas for $m\ge 2$ the dimension of $V_{m,n,\underline{c}}$ is
equal to the dimension of the space of homogeneous polynomials of
degree $n$ in $m-1$ variables and that is $\begin{pmatrix}m+n-2\\
n\end{pmatrix}$. In fact, if   $d(m,n)$ denotes the dimension of
the space of homogeneous polynomials of degree $n$ in $m$
variables we obtain the following recursion formula:
$$d(m,n)=\sum_{j=0}^n\ d(m-1,n-j).$$ 
Twice induction (first on
$n$, then on $m$) yields
 $$\sum_{j=0}^n\ \textstyle\begin{pmatrix}
                         m+j-2\\ j
                         \end{pmatrix}=
                    \textstyle     \begin{pmatrix}
                         m+n-1\\ n
                         \end{pmatrix}$$
and 
$$d(m,n)
=\sum_{j=0}^n\ d(m-1,j)
=\sum_{j=0}^n\ \textstyle \begin{pmatrix}
                                   m+j-2\\ j
                                   \end{pmatrix} \\
= \textstyle \begin{pmatrix}
                                   m+n-1\\ n
                     \end{pmatrix}.
$$
\qed

 The Ruelle zeta function $\zeta_{R}(\beta)$ for our special
choice of the parameter 
$\underline{\lambda}=\underline{\lambda}_{0}=(\frac{1}{2},\ldots,\frac{1}{2}$ 
has the form
 \[\zeta_{R}(\beta)
 =\prod_{k=0}^{m}\left(\det\left(1-\lambda_{0}^{k}\cL_{\beta}
                       \right)^{\binom{m}{k}}\right)^{(-1)^{k+1}}
 \]

as one easily checks. 
According to Proposition \ref{polydim} the contribution of the eigenspaces $\cB(D)_n$ 
for the eigenvalues $\varrho_{1,n}$ 
to the Ruelle function $\zeta_{R}(\beta)$  is given by
\begin{equation}\label{RuelleEFn}
\prod_{k=0}^{m}\prod_{r=0}^{\infty}\left(\left(1-\lambda_{0}^{k+r}\varrho_{1}\right)^{\binom{m}{k}
\binom{m+r-2}{r}}\right)^{(-1)^{k+1}}
\end{equation} 
with $k+r=n$.

\begin{Rem}
Note  that for any $\beta$ the  $\cB(D)_n$ represent the entire eigenspaces
for the eigenvalues $\varrho_{1,n}$. This is true for $\beta=0$ by 
Remark \ref{specL0} and Proposition \ref{polydim}. Indeed the eigenspace of the eigenvalue $\varrho = 2\lambda_{0}^{r}$ of the operator $\cL_{0}$ has dimension $\binom{m+r-1}{m-1}$ as one can easily check. On the other hand we have 
\[ \binom{m+r-1}{m-1}=\sum_{k=0}^{r}\binom{m+k-2}{k}\] where $\binom{m+k-2}{k}$ is just the dimension of the eigenspace of the eigenvalue $\varrho_{j,k}(\beta)=\varrho_{j}(\beta)\lambda_{0}^{k}$. But for $\beta = 0$ the eigenvalues $\varrho_{j}(\beta)$ are given by $2 \lambda_{0}^{j}$. 

Since the eigenvalue $\varrho_{1,k}(\beta)$ is holomorphic in the entire $\beta$-plane the dimension 
of its eigenspace does not depend on $\beta$ (see \cite{Kt66}, p.68) and is given by $\binom{m+k-2}{m-2}$.
\qed
\end{Rem}

We will check below that (\ref{RuelleEFn}) can be reduced  to the expression 
\begin{equation}\label{RuelleEFnsimple}
\frac{1-\lambda_{0}\varrho_{1}}{1-\varrho_{1}}.
\end{equation}
To see this one needs the following result on binomial coefficients.

\begin{Prop} \label{BINOM}For all $m\geqslant2$, all $r\geqslant 0$ and all $l \leqslant m-1$ the 
following identity holds
\[\sum_{k=0}^{\lfloor\frac{m}{2}\rfloor}\binom{m}{2k}\binom{m-2k+r}{l}=
\sum_{k=0}^{\lfloor\frac{m-1}{2}\rfloor}\binom{m}{2k+1}\binom{m-2k+r-1}{l},\]
where $\lfloor r\rfloor$ is the largest integer less or equal than $r$.
\end{Prop}

\pf
The proof is by induction on $m$, $r$ and $l$. For $m=2$, $r\geqslant 0$ and $l\in\{0,1\}$ one 
finds for the left respectively right hand side: $LHS=\binom{r+2}{l}+\binom{r}{l}$ respectively  
$RHS=2\binom{r+1}{l}$ and hence the two sides of the identity coincide for $l\in \{0,1\}$. 
Next we show that it suffices to show the identity for $r\ge 0$.
In fact, assume it holds for $r$. 
We show that then it holds 
for $r+1$ and hence for all $r$ . Since $\binom{n+1}{s}=\binom{n}{s}+\binom{n}{s-1}$ we find
\[\sum_{k=0}^{\lfloor\frac{m}{2}\rfloor}\binom{m}{2k}\binom{m-2k+r+1}{l}
=\sum_{k=0}^{[\frac{m}{2}]}\binom{m}{2k}\left(\binom{m-2k+r}{l}+\binom{m-2k+r}{l-1}\right).\]
 But the right hand side of this equation is equal to 
\begin{eqnarray*}
\sum_{k=0}^{\lfloor\frac{m-1}{2}\rfloor}\binom{m}{2k+1}\left(\binom{m-2k+r-1}{l}
        +\binom{m-2k+r-1}{l-1}\right)
=\sum_{k=0}^{\lfloor\frac{m-1}{2}\rfloor}\binom{m}{2k+1}\binom{m-2k+r}{l}
\end{eqnarray*}
which proves our claim.

Assume next the identity holds up to some $m$ and all $l\leqslant m-1$.
Then we show that the identity of Proposition \ref{BINOM} holds for $m+1$ and $0\leqslant l \leqslant m$.
 Indeed for the LHS of the identity one finds
\begin{eqnarray*}
&&\hskip -4em
\sum_{k=0}^{\lfloor\frac{m+1}{2}\rfloor}\binom{m+1}{2k}\binom{m-2k+1}{l}=\\
&=&
\sum_{k=0}^{\lfloor\frac{m+1}{2}\rfloor}\left(\binom{m}{2k}+\binom{m}{2k-1}\right)
           \binom{m-2k+1}{l}\\
&=&
\sum_{k=0}^{\lfloor\frac{m}{2}\rfloor}\binom{m}{2k}\binom{m-2k+1}{l}
    +\sum_{k=1}^{\lfloor\frac{m+1}{2}\rfloor}\binom{m}{2k-1}\binom{m-2k+1}{l}\\
&=&
\sum_{k=0}^{\lfloor\frac{m-1}{2}\rfloor}\binom{m}{2k+1}\binom{m-2k}{l}
+\sum_{k=0}^{\lfloor\frac{m-1}{2}\rfloor}\binom{m}{2k+1}\binom{m-2k-1}{l}\\
&=&
\sum_{k=0}^{\lfloor\frac{m-1}{2}\rfloor}\binom{m}{2k+1}\binom{m-2k}{l}
+\sum_{k=0}^{\lfloor\frac{m}{2}\rfloor}\binom{m}{2k}\binom{m-2k}{l}.
\end{eqnarray*}
But this is just
 \[\sum_{k=0}^{\lfloor\frac{m}{2}\rfloor}\binom{m+1}{2k+1}\binom{m-2k}{l} \] 
and hence the two sides of the identity coincide.
The proposition therefore holds for $m+1$, $0 \leqslant l \leqslant m-1$. 
We have still to show that it holds also for $l=m$, that means
\[\sum_{k=0}^{\lfloor\frac{m+1}{2}\rfloor}\binom{m+1}{2k}\binom{m+1-2k}{m}
=\sum_{k=0}^{\lfloor\frac{m}{2}\rfloor}\binom{m+1}{2k+1}\binom{m-2k}{m}.\]
But in this case we get $LHS = \binom{m+1}{m}$ respectively $RHS = m+1$ for this last equation, 
and hence the two sides agree. This proves the proposition.
\qed

The contribution (\ref{RuelleEFn}) of the eigenspace  $\cB(D)_n$ can be rewritten first as
\[\left(1-\lambda_{0}^{n}\varrho_{1}\right)^{\sum_{k=0}^{\lfloor\frac{m-1}{2}\rfloor}\binom{m}{2k+1}
\binom{m+n-2k-3}{n-2k-1}-\sum_{k=0}^{\lfloor\frac{m}{2}\rfloor}\binom{m}{2k}\binom{m+n-2k-2}{n-2k}}\]
or
\[\left(1-\lambda_{0}^{n}\varrho_{1}\right)^{\sum_{k=0}^{\lfloor\frac{m-1}{2}\rfloor}\binom{m}{2k+1}
\binom{m+n-2k-3}{m-2}-\sum_{k=0}^{\lfloor\frac{m}{2}\rfloor}\binom{m}{2k}\binom{m+n-2k-2}{m-2}}.\]
For the cases $n=0$ and $n=1$ one finds the factors $\frac{1}{1-\varrho_{1}}$ and 
$1-\lambda_{0}\varrho_{1}$. For $n \geqslant 2$, however, Proposition \ref{BINOM} shows that all the 
other eigenvalues $\varrho_{1,n}$ contribute the trivial factor $1$ to the zeta function. 
This proves 
$$\prod_{k=0}^{m}\prod_{r=0}^{\infty}\left(\left(1-\lambda_{0}^{k+r}\varrho_{1}\right)^{\binom{m}{k}
\binom{m+r-2}{r}}\right)^{(-1)^{k+1}}=
\frac{1-\lambda_{0}\varrho_{1}}{1-\varrho_{1}},
$$
i.e. the reduction of  (\ref{RuelleEFn}) to (\ref{RuelleEFnsimple}).
This factor leads to ``trivial'' zeros $\beta_{n}$ of
the Ruelle zeta function with $\beta_{n}= \frac{\log 2 + 2\pi i
n}{J}$ where $J=\sum_{l=1}^{m}J_{l}$. This conclusion can be drawn only if there are no 
cancellations with other eigenvalues, which at the moment we cannot exclude.

For the Kac-Baker model with this special
choice of parameters $\underline{\lambda}=\underline{\lambda}_{0}$ there exists another Ruelle transfer operator which in a certain sense is simpler than $\cL_{\beta}$. These parameters namely describe a Kac-Baker model with an interaction
consisting of one exponentially decreasing term only and whose
strength is just $J=\sum_{l=1}^{m}J_{l}$. The Ruelle transfer operator for this model is

\[\textstyle
\widetilde{\cL}_{\beta}g(z)= e^{\beta J z}
g(\frac{1}{2} + \frac{1}{2}\,z) + e^{-\beta J z}
g(-\frac{1}{2}+ \frac{1}{2} z), \] where $g$ is now holomorphic in a disc $D$ in the complex plane $\mathbb{C}$. 
The Ruelle zeta function in terms of the Fredholm determinants of
this operator has the form \[\zeta_{R}(\beta)=
\frac{\det(1-\frac{1}{2}\widetilde{\cL_{\beta}})}
{\det(1-\widetilde{\cL_{\beta}})}.\] 
It is easy to see that the
function $g_{1}(z)= \sinh\left( 2 \beta J z \right)$ is an
eigenfunction of the operator $\widetilde{\cL}_{\beta}$ with
eigenvalue $\varrho_{1}=\exp(\beta J)$. This eigenvalue
leads exactly to the zero of the Ruelle function we discussed
above.
 Indeed any eigenfunction
$f=f(\underline{z})$ of the Ruelle operator $\cL_{\beta}$
determines an eigenfunction $g=g(z)$ of the operator
$\widetilde{\cL}_{\beta}$ through $g(z):=f(z,\ldots,z)$ with the
same eigenvalue as long as this function does not vanish
identically. 
On the other hand given an eigenfunction $g_{k}=g_{k}(z)$ of the operator $\widetilde{\cL}_{\beta}$ with eigenvalue $\varrho_{k}$ we know that $g_{k}$ is holomorphic in the entire $\beta$-plane. 
Consider then the function $f_{k}=f_{k}(\underline{z})$ defined as
\[f_{k}(\underline{z}):= g_{k}\left(\frac{\sum_{l=1}^{m}J_{l}z_{l}}{J}\right)\] with $J=\sum_{l=1}^{m}J_{l}$. 
It is entire in $\mathbb{C}^{m}$ and since 
\[f_{k}\left(\pm\underline{\lambda}_{0}+\boldsymbol{\Lambda_{0}}\underline{z}\right)=g_{k}\left(\pm\lambda_{0}+\lambda_{0}\frac{\sum_{l=1}^{m}J_{l}z_{l}}{J}\right)\]
it is easy to see that the function $f_{k}$ is an eigenfunction of the operator $\cL_{\beta}$ with eigenvalue $\varrho_{k}$. 
But with $f_{k}$ also the functions $f_{k,n}(\underline{z}):= P_{n}(\underline{z})f_{k}(\underline{z}),n=0,1,...$ with $P_{n}$
 a polynomial homogeneous of degree $n$ in $m$ variables and invariant under translations $z_{l}\to z_{l}+c$ for $c\in \mathbb{R}$ are eigenfunctions of $\cL_{\beta}$ with eigenvalue $\varrho_{k,n}=\varrho_{k}\lambda_{0}^{n}$ as we have seen already for the eigenfunction $f_{1,n}$. Their degree of degeneracy is again given by $\binom{m+n-2}{n}$.
 The eigenfunctions $g_{1}(z)$ and $f_{1,n}(\underline{z})$ are only special cases of this quite general connection between the eigenfunctions and eigenvalues of the two operators $\cL_{\beta}$ and $\widetilde{\cL}_{\beta}$.
 An argument similar to the case $k=1$ shows that among all the eigenvalues $\varrho_{k,n}=\varrho_{k} \lambda_{0}^{n}$ only $\varrho_{k,n}$ with $n\in\{0,1\}$ give nontrivial contributions to the Ruelle zeta function, which coincide exactly with the contributions of the eigenvalues $\varrho_{k}$ of the operator $\widetilde{\cL}_{\beta}$ to the Ruelle zeta function for the model with $\underline{\lambda}=\underline{\lambda}_{0}$.
  
Our discussion shows that there exist infinitely many trivial zeros of
the Ruelle zeta function for the Kac model at least for this
special parameter $\underline{\lambda} = \underline{\lambda}_{0}$
as long as there are not an infinite number of cancellations
occuring among the eigenvalues, which one would certainly not expect. 
The structure of
the zeros of the dynamical zeta functions of this family of Kac-Baker models of
statistical mechanics hence seems to be very similar to the one
well known for arithmetic zeta functions like the Riemann
function. It would be interesting to determine their nontrivial zeros and poles and their distribution at least numerically.

\section{Matrix elements of the Kac-Gutzwiller operator $\widetilde{\mathcal{G}}_{\beta}$}\label{MatEl}

For the numerical determination of the zeros of the Ruelle zeta
function $\zeta_{R}(\beta)=\zeta_{R}(1,\beta)$ the modified
Kac-Gutzwiller operator $\widetilde{\mathcal{G}}_{\beta}$ seems to
be best suited. This was realized already in the case $m=1$ by
Gutzwiller in (see \cite[\S7]{Gu82}). For a numerical study of the
eigenvalues and eigenfunctions of a closely related Ruelle
operator see \cite{Thr94}. It turns out that also for general $m$
the matrix elements
$\widetilde{\mathbb{G}}_{\underline{\alpha},\underline{\delta}}$
of the operator $\widetilde{\mathcal{G}}_{\beta}$ can be
determined explicitly in the basis of
$L^{2}(\mathbb{R}^{m},d\underline{\xi})$ given by the Hermite
functions $\left\{h_{\underline{\alpha}}\right\}$. We start with
 the identity (see
\cite[(41)]{Gu82})
\begin{eqnarray*}
&&\hskip -2em \int_{\mathbb{R}^{m}}\int_{\mathbb{R}^{m}}
\prod_{i=1}^{m}\frac{1}{(4\pi\sinh\gamma_{i})^{\frac{1}{2}}}
\exp\left(-\frac{1}{4}\sum_{i=1}^{m}\left(\left(\xi_{i}^{2}+\eta_{i}^{2}\right)
\tanh\frac{\gamma_{i}}{2}+\frac{\left(\xi_{i}-\eta_{i}\right)^{2}}{\sinh\gamma_{i}}
\right)\right)\cdot\\
&&\cdot \exp\left(\pm\underline{R}\cdot\underline{\eta}-
\frac{\underline{\xi}^{2}}{4}-\frac{\underline{\rho}^{2}}{2}
+\underline{\rho}\cdot\underline{\xi}-\frac{\underline{\eta}^{2}}{4}
-\frac{\underline{\sigma}^{2}}{2}+\underline{\sigma} \cdot\underline{\eta}\right)
d\underline{\xi}d\underline{\eta}=\\
&&=
(2\pi)^{\frac{m}{2}}\exp\left(\frac{\underline{R}^{2}}{2}+\underline{\rho}
\cdot\boldsymbol{\Lambda}\underline{\sigma} \pm \underline{R}\cdot
\underline{\sigma} \pm \underline{R} \cdot \boldsymbol{\Lambda}
\underline{\rho}- \frac{1}{2} \sum_{i=1}^{m} \gamma_{i} \right),
\end{eqnarray*} where as before $ \boldsymbol{\Lambda}$ denotes the diagonal
matrix with entries $\boldsymbol{\Lambda}_{i,j}=
\lambda_{i}\delta_{i,j}$. Adding the two identities for $\pm
\underline{R}$ one gets
\begin{eqnarray*}
&&\hskip -2em 2 \int_{\mathbb{R}^{m}}\int_{\mathbb{R}^{m}}
\prod_{i=1}^{m}\frac{1}{(4\pi\sinh\gamma_{i})^{\frac{1}{2}}}
\exp\left(-\frac{1}{4}\sum_{i=1}^{m}\left(\left(\xi_{i}^{2}
+\eta_{i}^{2}\right)\tanh\frac{\gamma_{i}}{2}
+\frac{\left(\xi_{i}-\eta_{i}\right)^{2}}{\sinh\gamma_{i}}\right)\right)\cdot\\
&&\cdot  \cosh(\underline{R}\cdot\underline{\eta})\,
\exp\left(-\frac{\underline{\xi}^{2}}{4}-\frac{\underline{\rho}^{2}}{2}
+\underline{\rho}\cdot\underline{\xi}-\frac{\underline{\eta}^{2}}{4}
-\frac{\underline{\sigma}^{2}}{2}+\underline{\sigma} \cdot\underline{\eta}\right)
d\underline{\xi}\,d\underline{\eta}=\\
&&=2
\,(2\pi)^{\frac{m}{2}}\cosh\left(\underline{R}\cdot\underline{\sigma}
+\underline{R}\cdot\boldsymbol{\Lambda}\underline{\rho}\right)
\exp\left(\frac{\underline{R}^{2}}{2}
+\underline{\rho}\cdot\boldsymbol{\Lambda}\underline{\sigma} -
\frac{1}{2} \sum_{i=1}^{m} \gamma_{i} \right).
\end{eqnarray*}
But (see \cite[(37)]{Gu82})
\[ \exp\left(-\frac{\underline{\xi}^{2}}{4}
-\frac{\underline{\rho}^2}{2}-\underline{\rho}\cdot\underline{\xi}\right)
=(2\pi)^{\frac{m}{4}}\sum_{\underline{\alpha}\in
\mathbb{N}_{0}^{m}}
\frac{\underline{\rho}^{\underline{\alpha}}}{\sqrt{\underline{\alpha}!}}
\, h_{\underline{\alpha}}(\underline{\xi})\]
and hence
\begin{eqnarray*}
&&\hskip -2em \sum_{\underline{\alpha}\in
\mathbb{N}_{0}^{m}}\sum_{\underline{\delta}\in \mathbb{N}_{0}^{m}}
2 \prod_{i=1}^{m}\frac{1}{(4\pi\sinh\gamma_{i})^{\frac{1}{2}}}
\int_{\mathbb{R}^{m}}\int_{\mathbb{R}^{m}}
\exp\left(-\frac{1}{4}\sum_{i=1}^{m}\left(\left(\xi_{i}^{2}
+\eta_{i}^{2}\right)\tanh\frac{\gamma_{i}}{2}+\frac{\left(\xi_{i}
-\eta_{i}\right)^{2}}{\sinh\gamma_{i}}\right)\right)\cdot\\
&&\cdot\cosh(\underline{R}\cdot\underline{\eta})
\frac{\underline{\rho}^{\underline{\alpha}}}{\sqrt{\underline{\alpha}!}}
\frac{\underline{\sigma}^{\underline{\delta}}}{\sqrt{\underline{\delta}!}}
\,h_{\underline{\alpha}}(\underline{\xi})\,h_{\underline{\delta}}
(\underline{\eta})\,d\underline{\xi}\,d\underline{\eta}=\\
&&=2\, \exp\left(
\frac{\underline{R}^{2}}{2}+\underline{\rho}\cdot\boldsymbol{\Lambda}
\underline{\sigma}-\frac{1}{2}\sum_{i=1}^{m}\gamma_{i}\right)
\cosh\left(
\underline{R}\cdot\left(\underline{\sigma}+\boldsymbol{\Lambda}
\underline{\rho}\right)\right).
\end{eqnarray*}
Comparing this for $\underline{R}=\underline{R_{0}}:=\sqrt{\beta}
(\sqrt{J_1},\ldots,\sqrt{J_m})$ with the kernel
$\widetilde{\mathcal{G}}_{\beta}\left(\underline{\xi},\underline{\eta}\right)$
of the modified Kac-Gutzwiller operator in Theorem \ref{THM} we
hence find
\begin{eqnarray*}
&&\hskip-4em\sum_{\underline{\alpha}\in
\mathbb{N}_{0}^{m}}\sum_{\underline{\delta}\in
\mathbb{N}_{0}^{m}}\int_{\mathbb{R}^{m}}\int_{\mathbb{R}^{m}}
\widetilde{\mathcal{G}}_{\beta}\left(\underline{\xi},\underline{\eta}\right)
\,h_{\underline{\alpha}}(\underline{\xi})\,h_{\underline{\delta}}(\underline{\eta})
\frac{\underline{\rho}^{\underline{\alpha}}}{\sqrt{\underline{\alpha}!}}
\frac{\underline{\sigma}^{\underline{\delta}}}{\sqrt{\underline{\delta}!}}\,
d\underline{\xi}\, d\underline{\eta} =\\
&=& 2
\exp\left(\underline{\rho}\cdot\boldsymbol{\Lambda}\underline{\sigma}\right)
\cosh\left(\underline{R}_{0}\cdot
\left(\underline{\sigma}+\boldsymbol{\Lambda}\underline{\rho}\right)\right)
.
\end{eqnarray*}
In terms of the matrix elements
$\widetilde{\mathbb{G}}_{\underline{\alpha},\underline{\delta}}:=
\left(\widetilde{\mathcal{G}}_{\beta}h_{\underline{\alpha}},
h_{\underline{\delta}}\right)$
this reads
\begin{eqnarray*}
\sum_{\underline{\alpha}\in
\mathbb{N}_{0}^{m}}\sum_{\underline{\delta}\in
\mathbb{N}_{0}^{m}}\widetilde{\mathbb{G}}_{\underline{\alpha},\underline{\delta}}
\frac{\underline{\rho}^{\underline{\alpha}}}{\sqrt{\underline{\alpha}!}}
\frac{\underline{\sigma}^{\underline{\delta}}}{\sqrt{\underline{\delta}!}}
= 2
\exp\left(\underline{\rho}\cdot\boldsymbol{\Lambda}
\underline{\sigma}\right)\cosh\left(\underline{R}_{0}\cdot
\left(\underline{\sigma}+\boldsymbol{\Lambda}\underline{\rho}\right)\right)
.
\end{eqnarray*}
To determine from this the matrix elements
$\widetilde{\mathbb{G}}_{\underline{\alpha},\underline{\delta}}
=\widetilde{\mathbb{G}}_{\underline{\alpha},
\underline{\delta}}(\underline{R}_{0},\underline{\lambda})$
one first solves the problem
\begin{eqnarray*}
\sum_{\underline{\alpha}\in
\mathbb{N}_{0}^{m}}\sum_{\underline{\delta}\in \mathbb{N}_{0}^{m}}
\mathbb{A}_{\underline{\alpha},\underline{\delta}}
\frac{\underline{\rho}^{\underline{\alpha}}}{\sqrt{\underline{\alpha}!}}
\frac{\underline{\sigma}^{\underline{\delta}}}{\sqrt{\underline{\delta}!}}
= \exp\left(\underline{\rho}\cdot\boldsymbol{\Lambda}
\underline{\sigma}\right)
\exp\left(\underline{R}_{0}\cdot\left(\underline{\sigma}+\boldsymbol{\Lambda}
\underline{\rho}\right)\right)
.
\end {eqnarray*}
Obviously the right hand side factorizes into a product of exponentials
depending only on the $i$-th coordinates of the different variables.
Hence it suffices to solve the equation
\begin{eqnarray*}
\sum_{\alpha_{i}\in \mathbb{N}_{0}}\sum_{\beta_{i}\in
\mathbb{N}_{0}}
\mathbb{A}_{\alpha_{i},\beta_{i}}\frac{\rho_{i}^{\alpha_{i}}}{\sqrt{\alpha_{i}!}}
\frac{\sigma_{i}^{\beta_{i}}}{\sqrt{\beta_{i}!}} =
\exp\left(\rho_{i}\lambda_{i}\sigma_{i}\right)
\exp\left(R_{0,i}\left(\sigma_{i}+\lambda_{i}\rho_{i}\right)\right)
.
\end{eqnarray*}
Expanding the exponentials on the right hand side in $\sigma_{i}$ and $\rho_{i}$
leads to
\[
\mathbb{A}_{\alpha_{i},\beta_{i}}\left(R_{0,i},\lambda_{i}\right)
= \sqrt{\alpha_{i}!\beta_{i}!}
\,\lambda_{i}^{\alpha_{i}}\,R_{0,i}^{\mu_{i}}\sum_{k_{i}=1}^{M_{i}-\mu_{i}}
\frac{R_{0,i}^{2k_{i}}}{(M_{i}-\mu_{i}-k_{i})!\,k_{i}!\,(\mu_{i}+k_{i})!}
\]
where $M_{i}= \max\left\{\alpha_{i},\beta_{i}\right\}$ and
$\mu_{i}= \left |\alpha_{i}-\beta_{i}\right |$. It is not too
difficult to see that
 \[\sqrt{\alpha_{i}!\beta_{i}!}\sum_{k_{i}=0}^{M_{i}-\mu_{i}}
 \frac{R_{0,i}^{2k_{i}}}{(M_{i}-\mu_{i}-k_{i})!\,k_{i}!\,(\mu_{i}+k_{i})!}
 = \frac{1}{\sqrt{\alpha_{i}!\beta_{i}!}}\frac{M_{i}!}{\mu_{i}!}\,
 \Phi\left(\mu_{i}-M_{i},\mu_{i}+1;-R_{0,i}^{2}\right),
 \]
where $\Phi(\cdot,\cdot;\cdot)$ denotes the confluent
hypergeometric function. Therefore the matrix elements
$\mathbb{A}_{\underline{\alpha},\underline{\delta}} =\prod_{l=1}^m
\mathbb{A}_{\alpha_l,\delta_l}$ are given by
\[\mathbb{A}_{\underline{\alpha},\underline{\delta}} \left(\underline{R}_{0},
\underline{\lambda}\right) =
\frac{1}{\sqrt{\underline{\alpha}!\underline{\delta}!}}\,
\underline{\lambda}^{\underline{\alpha}}
\,\underline{R}_{0}^{\underline{\mu}}
\frac{\underline{M}!}{\underline{\mu}!}
\,\Phi\left(\underline{\mu}-\underline{M},
\underline{\mu}+\underline{1};-\underline{R}_{0}^{2}\right)
,
\]
where $\Phi\left(\underline{\mu}-\underline{M},\underline{\mu}
+\underline{1};-\underline{R}_{0}^{2}\right)
:= \prod_{i=1}^{m}\Phi\left(\mu_{i}-M_{i},\mu_{i}+1;-R_{0,i}^{2}\right).$

Finally we then find for the matrix elements
$\widetilde{\mathbb{G}}_{\underline{\alpha},\underline{\delta}}$
of the modified Kac-Gutzwiller operator
$\widetilde{\mathcal{G}}_{\beta}$
\[\widetilde{\mathbb{G}}_{\underline{\alpha},\underline{\delta}}
\left(\underline{R}_{0},\underline{\lambda}\right)=
\frac{1}{\sqrt{\underline{\alpha}!\underline{\delta}!}}
\,\underline{\lambda}^{\underline{\alpha}}\,\underline{R}_{0}^{\underline{\mu}}
\left(1+(-1)^{\left|\underline{\mu}\right|}\right)
\frac{\underline{M}!}{\underline{\mu}!}
\Phi\left(\underline{\mu}-\underline{M},\underline{\mu}
+\underline{1};-\underline{R}_{0}^{2}\right)
.
\]
This shows that
$\widetilde{\mathbb{G}}_{\underline{\alpha},\underline{\delta}}\neq
0$ only iff
$\left|\underline{\alpha}+\underline{\delta}\right|=0\mod 2$ which
generalizes the result for $m=1$ by Gutzwiller to arbitrary $m$.
 Since
$\left(\widetilde{\mathcal{G}}_{\beta}h_{\underline{\alpha}},
h_{\underline{\delta}}\right)=0$ if
$\left|\underline{\alpha}\right|\mod 2 \neq
\left|\underline{\delta}\right| \mod 2$
also for general $m$ the
Hilbert space $L^{2}(\mathbb{R}^{m},d\underline{\xi})$ can be
decomposed into two subspaces invariant under the operator
$\widetilde{\mathcal{G}}_{\beta}$ spanned by the Hermite functions
$\left\{h_{\underline{\alpha}}\right\}$ with
$\left|\underline{\alpha}\right|= 0\mod 2$ respectively
$\left|\underline{\alpha}\right|= 1\mod 2$. Obviously this
property of the operator $\widetilde{\mathcal{G}}_{\beta}$
corresponds to the fact that the Ruelle operator
$\mathcal{L}_{\beta}$ leaves invariant the subspaces of the Banach
space $\mathcal{B}(D)$ spanned by the functions
$F=F(\underline{z})$ which are even, respectively odd, under the
transformation $\underline{z}\to -\underline{z}$. This follows
from the fact that the Segal-Bargmann transform $B$ maps the
Hermite functions $h_{\underline{\alpha}}$ to the functions
$\zeta_{\underline{\alpha}}$ in (\ref{ZETA}) which under the above
transformation $\underline{z}\to-\underline{z}$ have parity
$\left|\underline{\alpha}\right|\mod 2$ as one checks easily. One
can use the representation of the operator
$\widetilde{\mathcal{G}}_{\beta}$ in terms of the matrix
$\widetilde{\mathbb{G}}_{\beta}$ to calculate its traces. For
instance one finds
\[\trace \widetilde{\mathbb{G}}_{\beta}=  \sum_{\underline{\alpha}\in \mathbb{N}_{0}^{m}}\widetilde{\mathbb{G}}_{\underline{\alpha},\underline{\alpha}}=
2\sum_{\underline{\alpha}\in \mathbb{N}_{0}^{m}}\underline{\lambda}^{\underline{\alpha}}\,\Phi(-\underline{\alpha},\underline{1};-\beta \underline{J}).\]
Because the confluent hypergeometric function  $\Phi(-n,1;x)$ is identical to the Laguerre polynomial $L_{n}(x)= L_{n}^{0}(x)$
 we get by using the generating function for these polynomials the result
 \[\trace\widetilde{\mathbb{G}}_{\beta}=\frac{2}{\prod_{i=1}^{m}(1-\lambda_{i})}
 \exp\left(\sum_{i=1}^{m}\frac{\beta J_{i}\lambda_{i}}{1-\lambda_{i}}\right)\]
 which coincides with $\trace \mathcal{L}_{\beta}$.

{\bf Acknowledgements}: This work has been supported by the Deutsche Forschungsgemeinschaft 
through the DFG Forschergruppe ``Zetafunktionen und lokal symmetrische Räume''.
 D.M. thanks the IHES in Bures sur Yvette for partial financial support and the kind hospitality 
extended to him during the preparation of the paper.

\end{document}